\documentclass[a4paper,10pt]{article}

\usepackage[OT2,T1]{fontenc}
\usepackage[french, english]{babel}

\usepackage{amsfonts}
\usepackage{amsmath,amssymb}
\usepackage{pstricks}
\usepackage{amsthm}
\usepackage[all]{xy}

\title{Le d\'efaut d'approximation forte dans les groupes lin\'eaires connexes}
\author{Cyril Demarche}
\date{\today}

\newtheorem{theo}{Th\'eor\`eme}[section]
\newtheorem*{thm}{Th\'eor\`eme}

\newtheorem{defi}[theo]{D\'efinition}
\newtheorem{prop}[theo]{Proposition}

\newtheorem{cor}[theo]{Corollaire}
\newtheorem{lem}[theo]{Lemme}
\newtheorem{rem}[theo]{Remarque}
\newenvironment{dem}{\noindent
  \textit{{D\'emonstration}} : }
  {\hfill \qedsymbol\newline}
\newenvironment{ex}{\noindent{\textbf{Exemples :}}}{}
\newenvironment{exs}{\noindent{\textbf{Exemples :}}}{}

\DeclareSymbolFont{rsfs}{U}{rsfs}{m}{n}
\DeclareSymbolFontAlphabet{\mathcal}{rsfs}

  \newcommand{\SSI}{si et seulement si }
  
\newcommand{\ensemble}[1]{\ensuremath{\mathbf #1} \xspace}
  \renewcommand{\H}{\ensemble H}
  \newcommand{\N}{\ensemble N}

  \newcommand{\Z}{\ensemble Z}
  
  \newcommand{\Q}{\ensemble Q}
  \newcommand{\R}{\ensemble R}
  \renewcommand{\P}{\ensemble P}

  \renewcommand{\L}{\ensemble L}
  \newcommand{\F}{\ensemble F}
  \renewcommand{\L}{\ensemble L}
  \newcommand{\G}{\ensemble G}
  \newcommand{\A}{\ensemble A}

  \newcommand\cyr[1]{{\fontencoding{OT2}\fontfamily{wncyr}\selectfont #1}}

\newenvironment{changemargin}[2]{
 \begin{list}{}{
  \setlength{\topsep}{0pt}
  \setlength{\leftmargin}{#1}
  \setlength{\rightmargin}{#2}
  \setlength{\listparindent}{\parindent}
  \setlength{\itemindent}{\parindent}
  \setlength{\parsep}{\parskip}
 }
\item[]}{\end{list}}

\newcommand{\fct}[4]{\begin{displaymath}
\begin{array}{ccc}
#1 & \longrightarrow & #2 \\
#3 & \longmapsto & #4
\end{array}
\end{displaymath} }

\begin{document}	

\selectlanguage{french}

\maketitle

\begin{abstract}
On consid\`ere un groupe alg\'ebrique lin\'eaire connexe $G$ d\'efini sur un corps de nombres $k$. On \'etablit une suite exacte d\'ecrivant l'adh\'erence du groupe $G(k)$ des points rationnels de $G$ dans le groupe des points ad\'eliques de $G$. On en d\'eduit que le d\'efaut d'approximation forte sur $G$ est mesur\'e par le groupe de Brauer alg\'ebrique de $G$ via l'obstruction de Brauer-Manin enti\`ere. On montre aussi que l'obstruction de Brauer-Manin enti\`ere sur un torseur sous $G$ est la seule obstruction \`a l'existence d'un point entier sur ce torseur. On obtient enfin une suite exacte (non-ab\'elienne) de type Poitou-Tate pour la cohomologie galoisienne du groupe $G$.
Les ingr\'edients principaux pour la preuve de ces r\'esultats sont les th\'eor\`emes de dualit\'e locale et globale pour les complexes de $k$-tores de longueur deux et les applications d'ab\'elianisation en cohomologie galoisienne d\'efinies par Borovoi.
\end{abstract}

\selectlanguage{english}

\begin{abstract}
Let $G$ be a connected linear algebraic group over a number field $k$. We establish an exact sequence describing the closure of the group $G(k)$ of rational points of $G$ in the group of adelic points of $G$. This exact sequence describes the defect of strong approximation on $G$ in terms of the algebraic Brauer group of $G$. In particular, we deduce from those results that the integral Brauer-Manin obstruction on a torsor under the group $G$ is the only obstruction to the existence of an integral point on this torsor. We also obtain a non-abelian Poitou-Tate exact sequence for the Galois cohomology of the linear group $G$. The main ingredients in the proof of those results are the local and global duality theorems for complexes of $k$-tori of length two and the abelianization maps in Galois cohomology introduced by Borovoi.
\end{abstract}

\selectlanguage{french}

\section{Introduction}

Soit $k$ un corps de nombres d'anneau des entiers $\mathcal{O}_k$, $\Omega_k$ l'ensemble des places de $k$, et $\overline{k}$ une cl\^oture alg\'ebrique de $k$. Si $v$ est une place de $k$, on note $k_v$ le compl\'et\'e de $k$ en $v$, et $\mathcal{O}_v$ l'anneau des entiers de $k_v$.
Par convention, on appelle $k$-groupe alg\'ebrique un sch\'ema en groupes de type fini sur $\textup{Spec } k$ (un tel sch\'ema en groupes est n\'ecessairement s\'epar\'e, voir \cite{SGA3}, expos\'e $\textup{VI}_{\textup{A}}$, section 0.2); pour un tel groupe $G$, on note indiff\'eremment $H^0(k, G)$ ou $G(k)$ l'ensemble des $k$-points de $G$.
Si $G / k$ est un groupe alg\'ebrique, on se donne un sch\'ema en groupes
s\'epar\'e $\mathcal{G}$ sur un ouvert $U$ (non vide) de
$\textup{Spec}(\mathcal{O}_k)$, tel que la fibre g\'en\'erique de $\mathcal{G}$ soit isomorphe \`a $G$, et on d\'efinit $P^0(k, G) :=
\prod'_{v \in \Omega_k} H^0(k_v, G)$ comme \'etant le produit
restreint des groupes $H^0(k_v, G)$ (avec la convention que si $v$ est
une place archim\'edienne, $H^0(k_v, G)$ d\'esigne le groupe des
composantes connexes de $G(k_v)$) par rapport aux sous-groupes
$H^0(\mathcal{O}_v, \mathcal{G})$. L'ensemble $P^0(k, G)$ ainsi obtenu
est ind\'ependant du mod\`ele $\mathcal{G}$ choisi. On munit ce groupe de sa topologie
de produit restreint, que l'on appelle indiff\'eremment
\emph{topologie ad\'elique} ou \emph{topologie forte}. On s'int\'eresse \`a \emph{l'adh\'erence forte} de $G(k)$ dans $P^0(k, G)$ : \'etant donn\'ee une partie $P$ de $P^0(k, G)$, on note $\overline{P}$ l'adh\'erence de $P$ dans $P^0(k, G)$ muni de cette topologie ad\'elique. On cherche notamment \`a d\'eterminer l'adh\'erence $\overline{G(k).\prod_{v \in S_0} H^0(k_v, G)}$ dans $P^0(k, G)$, o\`u $S_0$ est un ensemble fini de places de $k$ (on doit enlever un certain nombre de places de $k$, comme c'est le cas dans le th\'eor\`eme d'approximation forte pour les groupes semi-simples simplement connexes : voir \cite{PR}, th\'eor\`eme 7.12 ); ici, $G(k).\prod_{v \in S_0} H^0(k_v, G)$ d\'esigne exactement l'ensemble des produits d'un \'el\'ement de $G(k)$ par un \'el\'ement de $\prod_{v \in S_0} H^0(k_v, G)$ dans $P^0(k, G)$ (ou de fa\c con \'equivalente le sous-groupe engendr\'e par de tels produits : on montre en effet facilement que l'ensemble des produits est un sous-groupe de $P^0(k, G)$). 
On va d\'ecrire cette adh\'erence en termes du groupe de Brauer
cohomologique de $G$, d\'efini par
$\textup{Br}(G) := H^2_{\textup{\'et}}(G, \G_m)$; on d\'efinit aussi le
groupe de Brauer alg\'ebrique de $G$ par la formule $\textup{Br}_1(G) :=
\textup{Ker}(\textup{Br}(G) \rightarrow \textup{Br}(\overline{G}))$,
$\overline{G}$ d\'esignant le $\overline{k}$-groupe obtenu \`a partir
de $G$ par extension des scalaires de $k$ \`a $\overline{k}$. Enfin on
consid\`ere souvent le groupe quotient $\textup{Br}_a(G) :=
\textup{Br}_1(G) / \textup{Br}(k)$ (que l'on identifiera, via la section unit\'e de $G$, au sous-groupe de $\textup{Br}_1(G)$ form\'e des \'el\'ements nuls sur l'unit\'e de $G$), et ses sous-groupes finis
$\textup{\cyr{B}}(G) := \textup{Ker}(\textup{Br}_a(G) \rightarrow
\prod_{v \in \Omega_k} \textup{Br}_a(G \times_k k_v))$ et $\textup{\cyr{B}}_{\omega}(G) := \left\{ A \in \textup{Br}_a(G) : A_v
  = 0 \in \textup{Br}_a(G_v) \textup{ pour presque toute place } v \right\}$.
On utilise \'egalement les notations usuelles suivantes : si $H$ est un $k$-groupe alg\'ebrique (non n\'ecessairement commutatif), et $i = 0$ ou $1$, on note $H^i(k, H) := H^i(\Gamma_k, H(\overline{k}))$ l'ensemble de cohomologie galoisienne (voir \cite{Ser}, I.5.1), et on d\'efinit les ensembles suivants :
$$\textup{\cyr{SH}}^i(k, H) = \textup{\cyr{SH}}^i(H) := \textup{Ker}\left( H^i(k, H) \rightarrow \prod_{v \in \Omega_k} H^i(k_v, H) \right)$$
L'ensemble $\textup{\cyr{SH}}^1(k, H)$ est fini, et dans le cas o\`u le groupe $H$ est commutatif, le groupe $\textup{\cyr{SH}}^1(k, H)$ est le groupe de Tate-Shafarevich du groupe $H$.

Rappelons \'egalement que si $A$ est un groupe topologique ab\'elien, on note $A^D$ le groupe des morphismes de groupes continus $A \rightarrow \Q / \Z$. On munit ce groupe $A^D$ de la topologie compacte-ouverte.

On rappelle enfin qu'un $k$-groupe semi-simple est dit presque $k$-simple lorque ce groupe n'a pas de $k$-sous-groupe ferm\'e distingu\'e non trivial de dimension strictement positive.
Avec ces notations, l'un des r\'esultats principaux de ce texte est alors le suivant :
\begin{thm}[Th\'eor\`eme \ref{theo fin}]
Soit $G / k$ un groupe alg\'ebrique r\'eductif \footnote{par
  convention, dans tout ce texte, "r\'eductif" signifie "r\'eductif
  connexe", suivant la d\'efinition de \cite{SGA3}, Expos\'e XIX,
  1.6.}. Soit $S_0$ un ensemble fini de places de $k$. On note
$G^{\textup{sc}}$ le rev\^etement simplement connexe du groupe
d\'eriv\'e de $G$ et $\rho : G^{\textup{sc}} \rightarrow G$ le
morphisme naturel. On suppose que $(G^{\textup{sc}})^i_{S_0}$ est
non compact pour tout $k$-facteur presque $k$- simple
$(G^{\textup{sc}})^i$ de $G^{\textup{sc}}$. Il existe une bijection
fonctorielle $\textup{\cyr{B}}(G)^D \cong \textup{\cyr{SH}}^1(k, G)$
qui munit le second ensemble d'une structure de groupe
ab\'elien. Alors l'adh\'erence
$\overline{G(k).\rho(G^{\textup{sc}}_{S_0})}$ est un sous-groupe
distingu\'e de $P^0(k, G)$, et on a une suite exacte de groupes,
fonctorielle en $G$ :
$$1 \rightarrow \overline{G(k).\rho(G^{\textup{sc}}_{S_0})} \rightarrow P^0(k,G) \xrightarrow{\theta} (\textup{Br}_a G)^D \rightarrow \textup{\cyr{SH}}^1(k, G) \rightarrow 0$$
\end{thm}

Les fl\`eches apparaissant dans cette suite sont les morphismes suivants : la premi\`ere fl\`eche est l'injection \'evidente, la deuxi\`eme, not\'ee $\theta$, provient de l'accouplement de Brauer-Manin $P^0(k, G) \times \textup{Br}(G) \rightarrow \Q / \Z$ dont on rappelle la d\'efinition \`a la section \ref{subsection AF}, et la derni\`ere fl\`eche est la compos\'ee du morphisme dual de l'injection naturelle $\textup{\cyr{B}}(G) \rightarrow \textup{Br}_a G$ avec l'isomorphisme canonique de groupes ab\'eliens $\textup{\cyr{B}}(G)^D \cong \textup{\cyr{SH}}^1(k, G)$ (voir par exemple le th\'eor\`eme 8.5 de Sansuc dans \cite{San}).

L'un des ingr\'edients de la preuve de ce r\'esultat est le th\'eor\`eme d'approximation forte (voir th\'eor\`eme 7.12 de \cite{PR}, et d\'ebut de la section \ref{subsection AF}), d\^u notamment \`a Kneser et Platonov, et qui concerne les groupes semi-simples simplement connexes.

Ce r\'esultat g\'en\'eralise aux groupes r\'eductifs quelconques les r\'esultats r\'ecents suivants :
\begin{itemize}
\item Le th\'eor\`eme 4.5(b) de \cite{CTXF}, d\^u \`a Colliot-Th\'el\`ene et Xu Fei, qui traite le cas o\`u le groupe $G$ est semi-simple, et qui donne dans ce cas les trois premiers termes de la suite exacte du th\'eor\`eme.
  \item Le r\'esultat r\'ecent de Harari (th\'eor\`eme 2 de \cite{Har}), qui traite du d\'efaut d'approximation forte pour les $k$-tores (et plus g\'en\'eralement pour les $1$-motifs) : le lien entre ce th\'eor\`eme de Harari et celui que l'on se propose de d\'emontrer provient de l'isomorphisme classique $\textup{Br}_a(T) \cong H^2(k, \widehat{T})$ pour un $k$-tore $T$ (voir \cite{San}, lemme 6.9.(ii)).
\end{itemize}
En outre, la preuve propos\'ee ici n'est pas un ``d\'evissage'' \`a partir de ces r\'esultats : on les red\'emontre au passage. Dans la section \ref{section preuve 2}, on mentionne une autre preuve d'une partie du r\'esultat principal par d\'evissage et r\'eduction au cas des tores et au th\'eor\`eme 2 de \cite{Har}.

\begin{rem}
\label{rem esphom}
{\rm
\begin{itemize}
	\item Le groupe $\textup{Br}_a G$ apparaissant dans ce r\'esultat est en g\'en\'eral infini (c'est d\'ej\`a le cas pour un tore ou pour un groupe semi-simple non simplement connexe). Par cons\'equent, ce r\'esultat sous cette forme ne fournit pas un algorithme permettant de d\'ecider en un nombre fini d'\'etapes si un point ad\'elique est dans l'adh\'erence de $G(k).\rho(G^{\textup{sc}}_{S_0})$ ou non. Cependant, d'apr\`es la remarque due \`a Colliot-Th\'el\`ene \`a la fin de \cite{Har}, si l'ensemble $\mathcal{X}(\mathcal{O})$ est vide, seul un nombre fini d'\'el\'ements de $\textup{Br}_1(X)$ est n\'ecessaire pour obtenir une obstruction de Brauer-Manin enti\`ere. On peut alors se demander si l'on peut d\'eterminer de fa\c con effective un sous-groupe fini de $\textup{Br}_1(X)$ qui suffit. Des r\'esultats ont \'et\'e obtenus r\'ecemment par Xu Fei et Wei
Dasheng (voir \cite{XFDW}) dans cette direction, en ce qui concerne les groupes de type multiplicatif.
          \item Des g\'en\'eralisations des techniques utilis\'ees ici
            permettent d'\'etendre les th\'eor\`emes \'enonc\'es \`a des
            espaces homog\`enes de groupes connexes \`a stabilisateurs
            connexes ou ab\'eliens : le cas des espaces homog\`enes de groupes semi-simples simplement connexes est trait\'e par Colliot Th\'el\`ene et Xu Fei dans \cite{CTXF}, et on peut \'etendre ces r\'esultats \`a des espaces homog\`enes de groupes lin\'eaires connexes quelconques. On relie le d\'efaut
            d'approximation forte sur de tels espaces homog\`enes \`a
            un certain sous-groupe de leur groupe de Brauer, pouvant
            contenir des \'el\'ements transcendants (c'est-\`a-dire
            des \'el\'ements du groupe de Brauer qui ne deviennent pas
            triviaux quand on \'etend les scalaires \`a une cl\^oture alg\'ebrique). Ces g\'en\'eralisations seront trait\'ees dans un prochain texte.
\end{itemize}
}
\end{rem}

On d\'eduit de ce r\'esultat sur le d\'efaut d'approximation forte de nouvelles
estimations sur le nombre de classes d'un groupe r\'eductif (voir
section \ref{section nbre classes}).

Les preuves de ces r\'esultats reposent sur les th\'eor\`emes de
dualit\'e pour l'hypercohomologie des complexes de tores de longueur
$2$ obtenus dans \cite{Dem1}. On d\'eduit \'egalement de ces
th\'eor\`emes de dualit\'e une version de la suite de Poitou-Tate,
pour la cohomologie (non-ab\'elienne) d'un groupe lin\'eaire connexe,
qui prolonge la suite exacte de Kottwitz-Borovoi (voir \cite{BorAMS},
th\'eor\`eme 5.15). On rappelle que $\textup{UPic}(\overline{G})$ est d\'efini comme le complexe de modules galoisiens (en degr\'es $-1$ et $0$)
$$\textup{UPic}(\overline{G}) := [\overline{k}(G)^* / \overline{k}^* \rightarrow \textup{Div}(\overline{G})]$$
qui est quasi-isomorphe au dual du complexe de $k$-tores $[T_G^{\textup{sc}} \rightarrow T_G]$ par le corollaire 2.20 de \cite{BVH} ($T_G$ et $T_G^{\textup{sc}}$ sont des tores maximaux de $G^{\textup{red}}$ et $G^{\textup{sc}}$). Dans l'\'enonc\'e suivant, $H^2(k, G)$ d\'esigne l'ensemble de cohomologie non ab\'elienne, d\'efini dans \cite{FSS}, section 1 ou dans \cite{Bor2}, et $\sim$ est une certaine relation d'\'equivalence sur $H^2(k, G)$ :
\begin{thm}[Th\'eor\`eme \ref{PTNA}]
Soit $G / k$ un groupe lin\'eaire connexe.
\begin{itemize}
	\item Alors on a une suite exacte naturelle d'ensembles
          point\'es (la premi\`ere ligne est une suite exacte de groupes), fonctorielle en $G$ :
\begin{displaymath}
\xymatrix{
& (\textup{Br}_a G )^D \ar[d] & P^0(k, G) \ar[l] & P^0(k, G)^{\textup{Br}_a} \ar[l] & 0 \ar[l] \\
& H^1(k, G) \ar[r] & P^1(k, G) \ar[r] & (\textup{Pic } G )^D \ar[d] & \\
0 & (k[G]^* / k^*)^D \ar[l] & \bigoplus_{v \in \Omega_k} H^2(k_v, G) / \sim \ar[l] & H^2(k, G) / \sim \ar[l] &
}
\end{displaymath}
	\item Si de plus $G$ v\'erifie que $G^i_{S_0} := \prod_{v \in S_0} G^i(k_v)$  est non compact pour tout $k$-facteur presque $k$-simple $G^i$ de $G^{\textup{sc}}$, alors on peut identifier $ P^0(k, G)^{\textup{Br}_a}$ \`a $\overline{G^{\textup{scu}}_{S_0} G(k)}$.
	\item Dans tous les cas, on dispose de la suite exacte duale de groupes ab\'eliens :
\begin{displaymath}
\xymatrix{
0 \ar[r] & \widehat{G}(k)^{\wedge} \ar[r] & \left( \prod' \widehat{G}(k_v) \right)^{\wedge} \ar[r] & H^2_{\textup{ab}}(k, G)^D \ar[d] \\
& H^1_{\textup{ab}}(k, G)^D \ar[d] & \prod_{v \in \Omega_k} \textup{Pic}(G_v) \ar[l] & \textup{Pic}(G) \ar[l] \\
& \textup{Br}_a(G) \ar[r] & \left( \prod'_v \textup{Br}_a(G_v)
\right)_{\textup{tors}} \ar[r] & \left( H^0_{\textup{ab}}(k, G)^D \right)_{\textup{tors}} \ar[d] \\
0 & \textup{Ker}(\rho)(k)^D \ar[l] & \P^2(k,\textup{UPic}(\overline{G})) \ar[l] & \H^2(k, \textup{UPic}(\overline{G})) \ar[l]
}
\end{displaymath}
o\`u le produit restreint $\prod'_v \textup{Br}_a(G_v)$ est
consid\'er\'e par rapport aux sous-groupes $\H^1(\mathcal{O}_v,
[\widehat{\mathcal{T}_G} \rightarrow
\widehat{\mathcal{T}_G^{\textup{sc}}}])$.
\end{itemize}
\end{thm}

Le plan de ce texte est le suivant : on montre d'abord \`a la section
\ref{section AF complexes} un r\'esultat concernant l'approximation forte pour les
complexes de tores \`a l'aide des th\'eor\`emes de dualit\'e de
\cite{Dem1}.
Puis gr\^ace aux techniques d'ab\'elianisation de Borovoi, et \`a
l'aide du th\'eor\`eme d'approximation forte, on d\'ecrit le d\'efaut
d'approximation forte pour les groupes r\'eductifs dans le cas non
compact, en terme d'obstruction de Brauer-Manin (voir section
\ref{subsection AF}). Ensuite, on donne des applications de ces r\'esultats, avec
notamment des estimations sur le nombre de classes d'un groupe
r\'eductif (section \ref{section nbre classes}). Enfin, on d\'emontre \`a la section
\ref{section PTNA} l'existence d'une suite de Poitou-Tate non ab\'elienne pour un
groupe r\'eductif qui g\'en\'eralise la suite de Poitou-Tate pour les
tores et qui prolonge la suite exacte de Kottwitz-Borovoi (voir
\cite{BorAMS}, th\'eor\`eme 5.15). On termine ce texte par une
application concernant l'obstruction de Brauer-Manin enti\`ere sur un
espace principal homog\`ene sous un groupe alg\'ebrique connexe
(section \ref{section entier}).

\paragraph{Remerciements}
Je remercie tr\`es chaleureusement David Harari pour son aide et sa
patience. Je remercie \'egalement Mikhail Borovoi, Jean-Louis Colliot-Th\'el\`ene et Philippe Gille
pour leurs pr\'ecieux commentaires.

\section{Approximation forte pour les complexes de tores}
\label{section AF complexes}
Dans cette partie, \'etant donn\'e un complexe $C = [T_1
\xrightarrow{\rho} T_2]$ de $k$-tores, on d\'ecrit l'adh\'erence de
$\H^0(k, C)$ dans $\P^0(k,C) := \prod'_v \H^0(k_v,C)$ pour la
topologie de produit restreint (le produit restreint est pris par
rapport aux sous-groupes $\H^0(\mathcal{O}_v, \mathcal{C})$).

Dans toute cette section, $C = [T_1 \xrightarrow{\rho} T_2]$ est un complexe de tores sur $k$ qui
s'\'etend en un complexe de tores $\mathcal{C} = [\mathcal{T}_1
\xrightarrow{\rho} \mathcal{T}_2]$ sur $U = \textup{Spec } \mathcal{O}_{k,
  S}$, $S$ \'etant un ensemble fini de places de $k$. 
Dans toute la suite, on supposera \emph{le noyau $\textup{Ker}(\rho)$
  du morphisme de $k$-tores fini}. On renvoie \`a \cite{Dem1} pour les
g\'en\'eralit\'es sur l'hypercohomologie des complexes de tores.

\begin{lem}
\label{lem H^0_r}
Soit $v$ une place hors de $S$. On note $\H^0_r(k_v, C) := \H^0(k_v, C) / \H^0(\mathcal{O}_v, \mathcal{C})$.\\
Alors $\H^0_r(k_v, C)$ est un groupe discret, i.e. le sous-groupe $\H^0(\mathcal{O}_v, \mathcal{C})$ est ouvert dans $\H^0(k_v, C)$.
\end{lem}

\begin{dem}
On consid\`ere le diagramme commutatif suivant :
\begin{displaymath}
\xymatrix{
H^0(\mathcal{O}_v, \mathcal{T}_1) \ar[r] \ar[d] & H^0(\mathcal{O}_v, \mathcal{T}_2) \ar[r] \ar[d] & \H^0(\mathcal{O}_v, \mathcal{C}) \ar[r] \ar[d] & H^1(\mathcal{O}_v, \mathcal{T}_1) = 0 \ar[d] & \\
H^0(k_v, T_1) \ar[r] \ar[d] & H^0(k_v, T_2) \ar[r] \ar[d] & \H^0(k_v, C) \ar[r] \ar[d] & Q \ar[d] \ar[r] & 0 \\
H^0_r(k_v, T_1) \ar[r] \ar[d] & H^0_r(k_v, T_2) \ar[r] \ar[d] & \H^0_r(k_v, C) \ar[r] \ar[d] & Q \ar[d] & \\
0 & 0 & 0 & 0 &
}
\end{displaymath}
o\`u $Q$ est un sous-groupe fini (sous-groupe de $H^1(k_v, T_1)$). Le
groupe $H^1(\mathcal{O}_v, \mathcal{T}_1)$ est trivial par le
th\'eor\`eme de Lang, apr\`es restriction \`a la fibre sp\'eciale de $\textup{Spec}(\mathcal{O}_v)$. Les colonnes de ce diagramme sont exactes par d\'efinition, et les deux premi\`eres lignes sont exactes.
Une chasse au diagramme assure l'exactitude de la suite suivante :
$$H^0_r(k_v, T_1) \rightarrow H^0_r(k_v, T_2) \rightarrow \H^0_r(k_v, C) \rightarrow Q \rightarrow 0$$
Or les deux premiers groupes de cette suite sont discrets, et le groupe $Q$ est fini, donc cela assure que le groupe $\H^0_r(k_v, C)$ est discret.

\end{dem}

Notons $\P^0_S(C) := \prod'_{v \notin S} \H^0(k_v, C)$

\begin{lem}
\label{ind fini}
L'image de $\H^0(k, C)$ dans $\P^0_S(C) / {\prod_{v \notin S} \H^0(\mathcal{O}_v, C)} = \bigoplus_{v \notin S} \H^0_r(k_v, C)$ est d'indice fini.
\end{lem}

\begin{dem}
On sait d\'ej\`a que ce r\'esultat est vrai pour un tore (voir \cite{HarAF}, lemme 3).
On a un diagramme commutatif \`a lignes exactes (voir preuve du lemme \ref{lem H^0_r} pour la seconde ligne) de la forme :

\begin{displaymath}
\xymatrix{
H^0(k, T_2) \ar[r] \ar[d] & \H^0(k, C) \ar[r]^{h} \ar[d]^{f} & H^1(k, T_1) \ar[d]^{g} \ar[r] & H^1(k, T_2) \ar[d] \\
\bigoplus_{v \notin S} H^0_r(k_v, T_2) \ar[r] & \bigoplus_{v \notin S} \H^0_r(k_v, C) \ar[r]^{h'} & \bigoplus_{v \notin S} H^1(k_v, T_1) \ar[r] & \bigoplus_{v \notin S} H^1(k_v, T_2) 
}
\end{displaymath}
On sait que le conoyau du morphisme $g$ est fini : en effet, via la
suite exacte de Poitou-Tate pour $T_1$ (voir \cite{HSz}, th\'eor\`eme
5.6 ou \cite{Dem1}, th\'eor\`eme 6.1), on sait que le conoyau de $g' : H^1(k, T_1) \rightarrow \bigoplus_{v \in \Omega_k} H^1(k_v, T_1)$ s'identifie \`a un sous-groupe de $H^1(k, \widehat{T_1})^D$, donc $\textup{Coker}(g')$ est fini, et donc $\textup{Coker}(g)$ aussi puisque ce dernier est un quotient de $\textup{Coker}(g')$.

On s'int\'eresse au diagramme suivant, \`a lignes exactes :
\begin{displaymath}
\xymatrix{
0 \ar[r] & \textup{Ker } h \ar[r] \ar[d]^{f'} & \H^0(k, C) \ar[r] \ar[d]^{f} & \textup{Im } h \ar[r] \ar[d]^{\tilde{g}} & 0 \\
0 \ar[r] & \textup{Ker } h' \ar[r] & \bigoplus_{v \notin S} \H^0_r(k_v, C) \ar[r] & \textup{Im } h' \ar[r] & 0
}
\end{displaymath}

On va montrer que les conoyaux de $\tilde{g}$ et de $f'$ sont finis.
Pour $\tilde{g}$, cela r\'esulte de la finitude du noyau de $H^1(k,
T_2) \rightarrow \bigoplus_{v \notin S} H^1(k_v, T_2)$ (cette finitude
est une cons\'equence de la finitude de
$\textup{\cyr{SH}}^1_{\omega}(k, T_2)$) et de la finitude du conoyau
de $g$, par le lemme du serpent. Pour $f'$, c'est une cons\'equence
directe de la finitude du conoyau du morphisme $H^0(k, T_2)
\rightarrow \bigoplus_{v \notin S} H^0_r(k_v, T_2)$ (qui est une
cons\'equence de la finitude du nombre de classes de $T_2$ : voir \cite{PR}, th\'eor\`eme 5.1). On applique alors le lemme du serpent au diagramme pr\'ec\'edent pour montrer que $\textup{Coker } f$ est fini.
\end{dem}

\begin{prop}
\label{ramifie}
Pour un ouvert $U$ de $\textup{Spec}(\mathcal{O}_k)$ assez petit, on a une suite exacte (fonctorielle en $C$)
$$\H^0(U, \mathcal{C}) \rightarrow \H^0(k, C) \rightarrow \bigoplus_{v \notin S} \H^0_r(k_v, C)$$
\end{prop}

\begin{dem}
Cette suite est clairement un complexe. Pour l'exactitude, on s'int\'eresse au diagramme commutatif \`a lignes exactes suivant :
\begin{displaymath}
\label{diag tores nr}
\xymatrix{
H^0(U, \mathcal{T}_1) \ar[r] \ar[d] & H^0(U, \mathcal{T}_2) \ar[r] \ar[d] & H^0(U, \mathcal{C}) \ar[r] \ar[d] & H^1(U, \mathcal{T}_1) \ar[r] \ar[d] & H^1(U, \mathcal{T}_2) \ar[d] \\
H^0(k, T_1) \ar[r] \ar[d] & H^0(k, T_2) \ar[r] \ar[d] & \H^0(k, C) \ar[r] \ar[d] & H^1(k, T_1) \ar[d] \ar[r] & H^1(k, T_2) \\
\bigoplus_{v \notin S} H^0_r(k_v, T_1) \ar[r] & \bigoplus_{v \notin S} H^0_r(k_v, T_2) \ar[r] & \bigoplus_{v \notin S} \H^0_r(k_v, C) \ar[r] & \bigoplus_{v \notin S} H^1(k_v, T_1) &
}
\end{displaymath}
L'exactitude de la troisi\`eme ligne de ce diagramme a \'et\'e montr\'ee dans la preuve du lemme \ref{lem H^0_r}. On sait \'egalement que la deuxi\`eme colonne est exacte. Enfin, la quatri\`eme colonne de ce diagramme est aussi exacte par le th\'eor\`eme A.8 de \cite{GiP}.

Or, pour $U$ assez petit, le morphisme $ H^1(U, \mathcal{T}_2) \to H^1(k, T_2)$ est injectif (voir le corrigenda de \cite{HSz}). Une chasse au diagramme dans (\ref{diag tores nr}) assure alors que l'exactitude de la suite de la proposition (i.e. de la troisi\`eme colonne) est une cons\'equence de la surjectivit\'e du morphisme $H^0(k, T_1) \rightarrow \bigoplus_{v \notin S} H^0_r(k_v, T_1)$. Or cette surjectivit\'e est elle-m\^eme une cons\'equence de la finitude du nombre de classes d'un tore (voir th\'eor\`eme 5.1 de \cite{PR}), quitte \`a augmenter $S$, i.e. quitte \`a r\'eduire $U$.
\end{dem}

\begin{cor}
La suite suivante est exacte :
$$\H^0(U, \mathcal{C})^{\wedge} \rightarrow \H^0(k, C)^{\wedge} \rightarrow \left( \bigoplus_{v \notin S} \H^0_r(k_v, C) \right)^{\wedge}$$
\end{cor}

\begin{dem}
C'est une cons\'equence directe du lemme \ref{ind fini} et de la proposition \ref{ramifie} (voir preuve du corollaire 1 de \cite{HarAF}).
\end{dem}

Notons d\'esormais $\mathcal{P}^0_S(\mathcal{C}) := \prod_{v \in S} \H^0(k_v, C) \times \prod_{v \notin S} \H^0(\mathcal{O}_v, \mathcal{C})$.

\begin{prop}
\label{prop ouvert compl}
Pour $U$ assez petit, on a une suite exacte, fonctorielle en $\mathcal{C}$ :
$$\H^0(U, \mathcal{C})^{\wedge} \rightarrow \mathcal{P}^0_S(\mathcal{C})^{\wedge} \xrightarrow{\theta} \H^1(k, \widehat{C})^D$$
\end{prop}

\begin{dem}
On consid\`ere le diagramme suivant :
\begin{displaymath}
\xymatrix{
\H^0(U, \mathcal{C})^{\wedge} \ar[r] \ar[d] & \mathcal{P}^0_S(\mathcal{C})^{\wedge} \ar[r]^{\theta} \ar[d]^{i} & \H^1(k, \widehat{C})^D \ar[d]^{=} \\
\H^0(k, C)^{\wedge} \ar[r] \ar[d] & \P^0(C)^{\wedge} \ar[r]^{\theta} \ar[d] & \H^1(k, \widehat{C})^D \\
\left( \bigoplus_{v \notin S} \H^0_r(k_v, C) \right)^{\wedge} \ar[r]^{\cong} & \left( \bigoplus_{v \notin S} \H^0_r(k_v, C) \right)^{\wedge} & 
}
\end{displaymath}
La deuxi\`eme ligne de ce diagramme est exacte (extraite de la suite
de Poitou-Tate pour $C$ : voir la deuxi\`eme ligne de la suite exacte
du th\'eor\`eme 6.1 de \cite{Dem1}), ainsi que la premi\`ere colonne
(voir le corollaire pr\'ec\'edent). La deuxi\`eme colonne est
clairement un complexe. Il suffit donc de montrer l'injectivit\'e de
l'application $i : \mathcal{P}^0_S(\mathcal{C})^{\wedge} \rightarrow
\P^0(C)^{\wedge}$. Pour cela, gr\^ace \`a la preuve du lemme 5.21 de \cite{Dem1}, on sait que le morphisme $\P^0(C)_{\wedge} \rightarrow \P^0(C)^{\wedge}$ est injectif, donc il suffit d'avoir l'injectivit\'e de $\mathcal{P}^0_S(\mathcal{C})^{\wedge} \rightarrow \P^0(C)_{\wedge}$. 

Montrons l'injectivit\'e de $\mathcal{P}^0_S(\mathcal{C})_{\wedge}
\rightarrow \P^0(C)_{\wedge}$, \`a la mani\`ere du lemme 5.15 de
\cite{Dem1}, \`a la diff\'erence qu'ici il est n\'ecessaire d'avoir
recours à la cohomologie plate : la preuve du lemme 5.15 de \cite{Dem1} assure qu'il suffit, en consid\'erant le triangle exact suivant
$$\mathcal{C} \xrightarrow{n} \mathcal{C} \rightarrow \mathcal{C} \otimes^{\L} \Z / n \rightarrow \mathcal{C}[1]$$
de montrer que pour toute place $v$ hors de $S$ et pour tout entier $n$ non-nul, le morphisme 
$$\H^0_{\textup{fppf}}(\mathcal{O}_v, \mathcal{C} \otimes^{\L} \Z / n) \rightarrow \H^0_{\textup{fppf}}(k_v, C \otimes^{\L} \Z / n)$$
est injectif, puis d'utiliser le fait que $\mathcal{T}_1$ et
$\mathcal{T}_2$ sont lisses sur $U$ pour identifier les groupes
d'hypercohomologie $\H^0_{\textup{fppf}}(\mathcal{O}_v, \mathcal{C})$
et $\H^0_{\textup{\'et}}(\mathcal{O}_v, \mathcal{C})$ (voir
\cite{GBr}, III, th\'eor\`eme 11.7).

Soit $n \in \N^*$. On dispose du triangle exact suivant dans la
cat\'egorie d\'eriv\'ee associ\'ee \`a la cat\'egorie des complexes
born\'es de faisceaux fppf (voir \cite{Dem1}, lemme 2.3) :
$$_n \textup{Ker } \rho [2] \rightarrow \mathcal{C} \otimes^{\L} \Z / n \rightarrow T_{\Z/n}(\mathcal{C}) \rightarrow {_n\textup{Ker } \rho}[3]$$ 
On en d\'eduit donc le diagramme commutatif \`a lignes exactes suivant, pour toute place $v$ hors de $S$ :
\begin{displaymath}
\xymatrix{
H^2_{\textup{fppf}}(\mathcal{O}_v, { _n \textup{Ker } \rho}) \ar[r] \ar[d] & \H^0_{\textup{fppf}}(\mathcal{O}_v, \mathcal{C} \otimes^{\L} \Z / n) \ar[r] \ar[d] & \H^0_{\textup{fppf}}(\mathcal{O}_v, T_{\Z/n}(\mathcal{C})) \ar[d] \\
H^2_{\textup{fppf}}(k_v, { _n \textup{Ker } \rho}) \ar[r] & \H^0_{\textup{fppf}}(k_v, C \otimes^{\L} \Z / n) \ar[r] & \H^0_{\textup{fppf}}(k_v,T_{\Z/n}(C))
}
\end{displaymath}
Or le groupe $H^2_{\textup{fppf}}(\mathcal{O}_v, { _n \textup{Ker } \rho})$ est trivial par le lemme III.1.1.(a) de \cite{Mil}, puisque ${ _n \textup{Ker } \rho}$ est un sch\'ema en groupes fini plat sur $\textup{Spec}(\mathcal{O}_v)$. Donc l'injectivit\'e de la deuxi\`eme fl\`eche verticale est une cons\'equence de celle de la troisi\`eme. Pour montrer celle-ci, on utilise \`a nouveau le lemme III.1.1.(b) de \cite{Mil}, qui implique que le groupe $H^1_v(\mathcal{O}_v, T_{\Z/n}(\mathcal{C}))$ est trivial, car le sch\'ema en groupes $T_{\Z/n}(\mathcal{C})$ est fini plat sur $\mathcal{O}_v$. Donc on en d\'eduit que le morphisme 
$$\H^0_{\textup{fppf}}(\mathcal{O}_v, T_{\Z/n}(\mathcal{C})) \rightarrow \H^0_{\textup{fppf}}(k_v,T_{\Z/n}(C))$$
est injectif, ce qui implique bien, via le diagramme pr\'ec\'edent, que le morphisme 
$$\H^0_{\textup{fppf}}(\mathcal{O}_v, \mathcal{C} \otimes^{\L} \Z / n) \rightarrow \H^0_{\textup{fppf}}(k_v, C \otimes^{\L} \Z / n)$$
est injectif, pour tout $n$ et toute place $v \notin S$.
On a donc montr\'e (puisque la cohomologie fppf de $\mathcal{C}$ co\"incide avec sa cohomologie \'etale par lissit\'e) l'injectivit\'e du morphisme $\mathcal{P}^0_S(\mathcal{C})_{\wedge} \rightarrow \P^0(C)_{\wedge}$.

Pour conclure la preuve, il suffit d\'esormais de montrer que $\mathcal{P}^0_S(\mathcal{C})_{\wedge}$ est bien la compl\'etion profinie de $\mathcal{P}^0_S(\mathcal{C})$. Puisque la limite projective commute au produit (pour des groupes ab\'eliens), il suffit de montrer que $\H^0(k_v, C)_{\wedge} = \H^0(k_v, C)^{\wedge}$ et $\H^0(\mathcal{O}_v, \mathcal{C})_{\wedge} = \H^0(\mathcal{O}_v, \mathcal{C})^{\wedge}$. Mais ces r\'esultats sont clairs puisque $H^0(k_v, T_i)/n$, $H^0(\mathcal{O}_v, \mathcal{T}_i)/n$, $H^1(k_v, T_i)$ et $H^1(\mathcal{O}_v, \mathcal{T}_i)$ sont finis.
\end{dem}

Notons $\overline{\H^0(U, \mathcal{C})}$ l'adh\'erence de l'image de $\H^0(U, \mathcal{C})$ dans $\mathcal{P}^0_S(\mathcal{C})$.

\begin{theo}
\label{theo ouvert adh}
Pour $U$ suffisamment petit, on a une suite exacte, fonctorielle en $\mathcal{C}$ : 
$$0 \rightarrow \overline{\H^0(U, \mathcal{C})} \rightarrow \mathcal{P}^0_S(\mathcal{C}) \xrightarrow{\theta} \H^1(k, \widehat{C})^D$$
\end{theo}

\begin{dem}
Il est clair que cette suite est un complexe (voir th\'eor\`eme 6.1 de
\cite{Dem1}). On note $Q$ le quotient de $\mathcal{P}^0_S(\mathcal{C})$ par $\overline{\H^0(U, \mathcal{C})}$. C'est un espace topologique s\'epar\'e. Montrons que le morphisme de compl\'etion $Q \rightarrow Q^{\wedge}$ est injectif. Pour cela, il suffit que $Q$ soit compactement engendr\'e. Or si $v \in S$, la finitude de $H^1(k_v, T_1)$ assure que $\H^0(k_v, C)$ contient un sous-groupe d'indice fini, qui est un quotient topologique de $T_2(k_v)$. Or $T_2(k_v)$ est engendr\'e par une partie compacte, donc $Q$ aussi. D'o\`u le diagramme commutatif
\begin{displaymath}
\xymatrix{
0 \ar[r] & \overline{\H^0(U, \mathcal{C})} \ar[r] \ar[d] & \mathcal{P}^0_S(\mathcal{C}) \ar[r] \ar[d] & Q \ar[d] \\
& \overline{\H^0(U, \mathcal{C})}^{\wedge} \ar[r] & \mathcal{P}^0_S(\mathcal{C})^{\wedge} \ar[r] & Q^{\wedge}
}
\end{displaymath}
Alors une chasse au diagramme, l'application des r\'esultats pr\'ec\'edents et l'injectivit\'e de la derni\`ere fl\`eche verticale assurent le r\'esultat.
\end{dem}

\begin{lem}
\label{lem meme image}
$\P^0(C)$ et $\P^0(C)^{\wedge}$ ont m\^eme image dans $\H^1(k, \widehat{C})^D$.
\end{lem}

\begin{dem}
Soit $\mathcal{C} := \left[ \mathcal{T}_1 \rightarrow \mathcal{T}_2 \right]$ un complexe de sch\'emas en groupes commutatifs plats de type fini \'etendant $C$ sur $\textup{Spec } \mathcal{O}_k$. On d\'efinit le groupe des classes de $\mathcal{C}$ par $\textup{Cl}(\mathcal{C}) := \P^0(k, C) / (\H^0(k, C).\mathcal{C}(\A(\infty)))$, o\`u $\mathcal{C}(\A(\infty)) := \prod_{v \in S_{\infty}} \H^0(k_v, C) \times \prod_{v \notin S_{\infty}} \H^0(\mathcal{O}_v, \mathcal{C})$.

\begin{lem}
Le groupe $\textup{Cl}(\mathcal{C})$ est fini.
\end{lem}

\begin{dem}
On consid\`ere le diagramme commutatif \`a lignes exactes suivant :
\begin{displaymath}
\xymatrix{
H^0(k, T_2) \times \mathcal{T}_2(\A({\infty})) \ar[r] \ar[d]^{a} & \H^0(k, C) \times \mathcal{C}(\A({\infty})) \ar[r] \ar[d]^{b} & H^1(k, T_1) \times \mathcal{Q} \ar[r] \ar[d]^{c} & H^1(k, T_2) \ar[d]^{d} \\
P^0(k, T_2) \ar[r] & \P^0(k, C) \ar[r] & P^1(k, T_1) \ar[r] & P^1(k, T_2)
}
\end{displaymath}
o\`u $\mathcal{Q}$ est le quotient de $\mathcal{C}(\A({\infty}))$ par l'image de $\mathcal{T}_2(\A({\infty}))$.
Le lemme du serpent assure que la finitude de $\textup{Coker}(b)$ est une cons\'equence des finitudes de $\textup{Coker}(a)$, $\textup{Coker}(c)$ et $\textup{Ker}(d)$. Or $\textup{Coker}(a)$ est fini par finitude du nombre de classes du tore $T_2$, $\textup{Coker}(c)$ est fini gr\^ace \`a la suite de Poitou-Tate pour $T_1$ (finitude de $H^1(k, \widehat{T_1})$), et $\textup{Ker}(d) = \textup{\cyr{SH}}^1(T_2)$ est fini. D'o\`u le lemme.
\end{dem}

Montrons que ce lemme implique le lemme \ref{lem meme image}. Il est clair qu'il suffit de montrer que le quotient $\P^0(C) / \overline{\H^0(k, C)}$ est compact. Or $\mathcal{C}(\A({\infty}))$ est un sous-groupe ouvert de $\P^0(k, C)$, donc $\overline{H^0(k, C)} \subset \H^0(k, C).\mathcal{C}(\A({\infty}))$, et on a une suite exacte de groupes ab\'eliens :
$$\mathcal{C}(\A({\infty})) \rightarrow \P^0(C) / \overline{\H^0(k, C)} \rightarrow \textup{Cl}(\mathcal{C}) \rightarrow 0$$
On remarque que le groupe $\mathcal{C}(\A({\infty}))$ est compact (par le th\'eor\`eme de Tychonov et la convention pour les groupes de cohomologie modifi\'es aux places infinies), et le groupe $\textup{Cl}(\mathcal{C})$ est fini par le lemme pr\'ec\'edent, donc la suite exacte assure que le groupe $\P^0(C) / \overline{\H^0(k, C)}$ est compact. Cela conclut la preuve du lemme \ref{lem meme image}.

\end{dem}

\begin{theo}
\label{theo AF complexe}
On note $\overline{\H^0(k, C)}$ l'adh\'erence forte de $\H^0(k, C)$
dans $\P^0(C)$. Alors on a une suite exacte, fonctorielle en $C$ :
$$0 \rightarrow \overline{\H^0(k, C)} \rightarrow \P^0(C) \rightarrow \H^1(k, \widehat{C})^D \rightarrow \textup{\cyr{SH}}^1(C) \rightarrow 0$$
\end{theo}

\begin{dem}
On choisit d'abord $U$ suffisamment petit, et $S$ correspondant. Le d\'ebut de la suite exacte (trois premiers termes) provient du th\'eor\`eme \ref{theo ouvert adh}, en passant \`a la limite inductive sur $T$ fini contenant $S$. On extrait la suite exacte suivante de la suite de Poitou-Tate :
$$\P^0(C)^{\wedge} \rightarrow  \H^1(k, \widehat{C})^D \rightarrow \textup{\cyr{SH}}^1(C) \rightarrow 0$$
et on conclut avec le lemme \ref{lem meme image}.
\end{dem}

\section{Approximation forte dans les groupes lin\'eaires connexes}
L'objectif de cette section est de d\'eduire le th\'eor\`eme principal
(th\'eor\`eme \ref{theo fin}) du th\'eor\`eme d'approximation forte et
du th\'eor\`eme \ref{theo AF complexe}. On commence par traiter les
places r\'eelles.

\subsection{Hypercohomologie modifi\'ee \`a la Tate et ab\'elianisation}
\label{section Tate}
L'objectif de cette section est d'\'etendre les r\'esultats connus, d\^us \`a Deligne dans \cite{Del}, \`a Breen dans \cite{Bre} et \`a Borovoi dans \cite{BorAMS}, sur l'ab\'elianisation de la cohomologie galoisienne (en degr\'e $0$ essentiellement) \`a la cohomologie modifi\'ee "\`a la Tate". 
Soient $G_1$ et $G_2$ deux groupes alg\'ebriques sur $\R$, $\varphi : G_1 \rightarrow G_2$ un morphisme de $\R$-groupes.

\begin{defi}
On note $Z^0(\R, [G_1 \rightarrow G_2])$ la cat\'egorie des couples $(X, P)$, o\`u $X$ est un $\R$-torseur sous $G_1$ et $P$ est une section du $\R$-torseur $X'$ obtenu en poussant $X$ par le morphisme $\varphi : G_1 \rightarrow G_2$, les fl\`eches \'etant les suivantes : deux objets $(X, P)$ et $(Y, Q)$ sont isomorphes s'il existe un isomorphisme de torseurs $\psi : X \rightarrow Y$ tel que l'image de $P$ par l'isomorphisme $\psi' : X' \rightarrow Y'$ (induit par $\psi$) est dans la m\^eme composante connexe de $Y'(\R)$ que $Q$.
L'ensemble des classes d'isomorphisme ainsi obtenues est not\'e $\widehat{\H}^0(\R, [G_1 \rightarrow G_2])$. C'est un ensemble point\'e par la classe du torseur trivial sous $G_1$ muni de la trivialisation fournie par le neutre de $G_2$.
\end{defi}

\begin{rem}
{\rm
On dispose de fl\`eches canoniques surjectives \'evidentes $\H^0(\R,
[G_1 \rightarrow G_2]) \rightarrow \widehat{\H}^0(\R, [G_1 \rightarrow
G_2])$, o\`u $\H^0(\R, [G_1 \rightarrow G_2])$ est d\'efini par
exemple dans \cite{BorAMS}, section 3.1.1.
}
\end{rem}

\begin{rem}
\label{rem cocycles}
{\rm
On dispose \'egalement d'une d\'efinition de l'ensemble
$\widehat{\H}^0(\R, [G_1 \rightarrow G_2])$ en terme de cocycles : on
note $Z^0$ l'ensemble des couples $(z, g)$ o\`u $z : \Gamma_R
\rightarrow G_1$ est un $1$-cocycle au sens usuel, et $g_2 \in G_2$,
de sorte que l'on ait la relation $\varphi(z_{\sigma}^{-1}) g_2 =
^{\sigma} g_2$ dans $G_2$. On munit alors $Z^0$ de la relation d'\'equivalence suivante : $(z, g_2) \simeq (z', g_2')$ s'il existe $g_1 \in G_1$ tel que $z'_{\sigma} = g_1^{-1} z_{\sigma} ^{\sigma} g_1$ et $g_2' g_2^{-1} \varphi(g_1) \in G_2^0(\R)$. Alors on d\'efinit $\widehat{\H}^0(\R, [G_1 \rightarrow G_2])$ comme \'etant le quotient de $Z^0$ par cette relation d'\'equivalence. On v\'erifie ais\'ement que les deux ensembles de cohomologie modifi\'es ainsi d\'efinis sont canoniquement isomorphes.
}
\end{rem}

\begin{rem}
{\rm L'ensemble $\widehat{\H}^0(\R, [T^{\textup{sc}} \rightarrow T])$ d\'efini  dans cette section s'identifie canoniquement au groupe d'hypercohomologie modifi\'ee du complexe de tores $[T^{\textup{sc}} \rightarrow T]$ not\'e \'egalement $\widehat{\H}^0(\R, [T^{\textup{sc}} \rightarrow T])$. Cette identification est \'evidente via la description en termes de cocycles de la remarque \ref{rem cocycles} qui co\"incide avec la d\'efinition usuelle utilisant des r\'esolutions compl\`etes du groupe $\Gamma_R$.}
\end{rem}

\begin{lem}
On dispose d'une suite exacte naturelle
$$\widehat{H}^0(\R, G_1) \rightarrow \widehat{H}^0(\R, G_2) \rightarrow \widehat{\H}^0(\R, [G_1 \rightarrow G_2]) \rightarrow H^1(\R, G_1) \rightarrow H^1(\R, G_2)$$
\end{lem}

\begin{dem}
La d\'efinition des fl\`eches apparaissant dans cette suite est claire : la premi\`ere est la fl\`eche usuelle induite par $\varphi$; la deuxi\`eme associe \`a une composante connexe de $G_2(\R)$ la classe du torseur trivial sous $G_1$ muni de la trivialisation de $G_2$ donn\'ee par un point de cette composante connexe; la troisi\`eme fl\`eche est juste la fl\`eche d'oubli de la trivialisation; la derni\`ere fl\`eche est la fl\`eche \'evidente induite par $\varphi$. Par d\'efinition des fl\`eches, la suite du lemme est bien un complexe. 

Montrons son exactitude :
\begin{itemize}
	\item Exactitude en $\widehat{H}^0(\R, G_2)$ :
Soit $P \in G_2(\R)$ dont la classe dans $\widehat{H}^0(\R, G_2)$ est d'image triviale dans $\widehat{\H}^0(\R, [G_1 \rightarrow G_2])$. Alors il existe un isomorphisme de $G_1$-torseurs $\psi : G_1 \rightarrow G_1$ induisant un morphisme $\psi' : G_2 \rightarrow G_2$ compatible avec $\varphi$, de sorte que $\psi'(P)$ est dans la composante neutre de $G_2(\R)$. Alors, si $e_1$ d\'esigne le neutre de $G_1$, $\psi^{-1}(e_1)$ est un point de $G_1(\R)$ dont l'image par $\varphi$ est dans la m\^eme composante connexe de $G_2(\R)$ que $P$. Donc la classe de $P$ dans $\widehat{H}^0(\R, G_2)$ se rel\`eve dans $\widehat{H}^0(\R, G_1)$.
	\item Exactitude en $\widehat{\H}^0(\R, [G_1 \rightarrow G_2])$ :
Soit $[(X,P)]$ une classe dans $\widehat{\H}^0(\R, [G_1 \rightarrow G_2])$ d'image triviale dans $H^1(\R, G_1)$. Cela signifie exactement que le torseur $X$ est isomorphe au torseur trivial sous $G_1$, et donc que la classe $[(X, P)]$ est \'egale \`a une classe $[(G_1, P')]$, $P' \in G_2(\R)$ correspondant \`a la trivialisation $P$ via une trivialisation $X \cong G_1$, donc $[(X, P)]$ est bien l'image de la classe de $P'$ dans $\widehat{H}^0(\R, G_2)$.
	\item Exactitude en $H^1(\R, G_1)$ : \'evident par d\'efinition de $\widehat{\H}^0(\R, [G_1 \rightarrow G_2])$.
\end{itemize}

\end{dem}

\begin{lem}
Soit 
$$1 \rightarrow [A \rightarrow B] \rightarrow [C \rightarrow D] \rightarrow [E \rightarrow F] \rightarrow 1$$
une suite exacte de modules crois\'es de $\R$-groupes. Alors on a une suite exacte
$$\widehat{\H}^0(\R, [A \rightarrow B]) \rightarrow \widehat{\H}^0(\R, [C \rightarrow D]) \rightarrow \widehat{\H}^0(\R, [E \rightarrow F]) \rightarrow \H^1(\R, [A \rightarrow B])$$
\end{lem}

\begin{dem}
L'exactitude en $\widehat{\H}^0(\R, [E \rightarrow F])$ se d\'eduit imm\'ediatement du r\'esultat analogue pour les groupes non modif\'es $\H^0(\R, .)$. Reste donc \`a montrer l'exactitude en $\widehat{\H}^0(\R, [C \rightarrow D])$. Soit $[(X, P)] \in \widehat{\H}^0(\R, [C \rightarrow D])$ d'image triviale dans $\widehat{\H}^0(\R, [E \rightarrow F])$.
On utilise alors les deux suites exactes suivantes 
$$H^1(\R, A) \rightarrow H^1(\R, C) \rightarrow \H^1(\R, E)$$
et
$$\widehat{H}^0(\R, B) \rightarrow \widehat{H}^0(\R, D) \rightarrow \widehat{H}^0(\R, F)$$
pour montrer que la classe du torseur $X$ provient d'une classe dans $H^1(\R, A)$, et quitte \`a remplacer $P$ par un autre point dans sa composante connexe, pour montrer que $P$ se rel\`eve dans $B(\R)$. Cela assure que la classe $[(X, P)]$ provient d'une classe dans $\widehat{\H}^0(\R, [A \rightarrow B])$.
\end{dem}

Soit d\'esormais $G$ un $\R$-groupe r\'eductif. On d\'efinit $\rho : G^{\textup{sc}} \rightarrow G$, $T$ et $T^{\textup{sc}}$ comme dans l'introduction de \cite{BorAMS}.

\begin{lem}
\begin{itemize}
	\item Les modules crois\'es $[G^{\textup{sc}} \rightarrow G]$ et $[T^{\textup{sc}} \rightarrow T]$ sont quasi-isomorphes.
	\item Le morphisme $\widehat{\H}^0(\R, [T^{\textup{sc}} \rightarrow T]) \rightarrow \widehat{\H}^0(\R, [G^{\textup{sc}} \rightarrow G])$ est une bijection d'ensembles point\'es.
\end{itemize}
\end{lem}

\begin{dem}
\begin{itemize}
	\item Voir \cite{BorAMS}, lemme 3.8.1.
	\item On dispose de deux suites exactes de complexes de groupes (en utilisant le premier point) :
$$1 \rightarrow \textup{Ker}(\rho)[1] \rightarrow [G^{\textup{sc}} \rightarrow G] \rightarrow \textup{Coker}(\rho) \rightarrow 1$$
et
$$1 \rightarrow \textup{Ker}(\rho)[1] \rightarrow [T^{\textup{sc}} \rightarrow T] \rightarrow \textup{Coker}(\rho) \rightarrow 1$$
On v\'erifie d'abord que les \'el\'ements du noyau de
$\widehat{H}^0(\R, \textup{Ker}(\rho)[1]) \rightarrow \widehat{\H}^0(\R,
[G^{\textup{sc}} \rightarrow G])$ sont exactement ceux qui se rel\`event
dans $\widehat{H}^{-1}(\R, \textup{Coker}(\rho))$ : en effet, on dispose
du diagramme commutatif exact naturel suivant :
\begin{displaymath}
\xymatrix{
& \widehat{H}^{-1}(\R, \textup{Coker } \rho) \ar[d] & & \\
\widehat{H}^0(\R, G^{\textup{sc}}) \ar[r] \ar[d]^{=} &
\widehat{H}^0(\R, \textup{Im } \rho) \ar[r] \ar[d] & H^1(\R,
\textup{Ker } \rho) \ar[r] \ar[d] & H^1(\R, G^{\textup{sc}})
\ar[d]^{=} \\
\widehat{H}^0(\R, G^{\textup{sc}}) \ar[r] &
\widehat{H}^0(\R, G) \ar[r] & \widehat{H}^0(\R,
[G^{\textup{sc}} \rightarrow G]) \ar[r] & H^1(\R, G^{\textup{sc}})
}
\end{displaymath}
et on conclut par une chasse au diagramme. On applique alors le lemme
pr\'ec\'edent, et on obtient le diagramme suivant, \`a lignes exactes :
\begin{changemargin}{-1cm}{1cm}
\begin{displaymath}
\xymatrix{
\widehat{H}^{-1}(\R, \textup{Coker}(\rho)) \ar[r] \ar[d]^{=} & H^1(\R, \textup{Ker}(\rho)) \ar[r] \ar[d]^{=} & \widehat{\H}^0(\R, [T^{\textup{sc}} \rightarrow T]) \ar[r] \ar[d] & \widehat{H}^0(\R, \textup{Coker}(\rho)) \ar[r] \ar[d]^{=} & H^2(\R, \textup{Ker}(\rho)) \ar [d]^{=} \\
\widehat{H}^{-1}(\R, \textup{Coker}(\rho)) \ar[r] & H^1(\R, \textup{Ker}(\rho)) \ar[r] & \widehat{\H}^0(\R, [G^{\textup{sc}} \rightarrow G]) \ar[r] & \widehat{H}^0(\R, \textup{Coker}(\rho)) \ar[r] & H^2(\R, \textup{Ker}(\rho))
}
\end{displaymath}
\end{changemargin}
et on conclut gr\^ace au lemme des cinq.
\end{itemize}
\end{dem}

\subsection{Approximation forte dans les groupes lin\'eaires connexes}
\label{subsection AF}
Commen\c cons par une d\'efinition qui nous servira dans toute la
suite :
\begin{defi}
Soit $G$ un $k$-groupe alg\'ebrique semi-simple simplement connexe, et
$S_0$ un ensemble fini de places de $k$. On
dit que le groupe $G$ v\'erifie l'approximation forte hors de $S_0$
lorsque $G(k).G_{S_0}$ est dense dans $P^0(k,G)$ (muni de la topologie ad\'elique).
\end{defi}

On rappelle ici l'\'enonc\'e du th\'eor\`eme d'approximation forte (voir \cite{PR}, th\'eor\`eme 7.12) 

\begin{thm}[Kneser, Platonov]
Soit $k$ un corps de nombres, $G / k$ un groupe alg\'ebrique semi-simple, simplement connexe. Soit $S_0$ un ensemble fini de places tel que $G^i_{S_0} := \prod_{v \in S_0} G^i(k_v)$  est non compact pour tout $k$-facteur presque $k$-simple $G^i$ de $G$. Alors $G$ v\'erifie l'approximation forte hors de $S_0$.
\end{thm}

En particulier, si le groupe $G$ est presque $k$-simple simplement connexe, et si $v_0$ est une place de $k$ telle que $G(k_{v_0})$ est non compact, alors $G(k).G(k_{v_0})$ est dense dans $P^0(k,G)$.

Desormais, dans cette section, $G$ est un groupe alg\'ebrique sur $k$, suppos\'e reductif. On note $G^{\textup{ss}}$ son sous-groupe d\'eriv\'e, $G^{\textup{sc}}$ le rev\^etement semi-simple simplement connexe de $G^{\textup{ss}}$, $T$ un tore maximal de $G$, $T^{\textup{sc}}$ l'image r\'eciproque de $T$ dans $G^{\textup{sc}}$ et $C := \left[ T^{\textup{sc}} \xrightarrow{\rho} T \right]$ le complexe de tores associ\'e. On se donne $S_1$ un ensemble fini de places, contenant les places archim\'ediennes de $k$, tel que $G$ s'\'etende respectivement en un sch\'ema en groupes r\'eductif $\mathcal{G}$ sur $V := \textup{Spec } \mathcal{O}_{k, S_1}$ (au sens de \cite{SGA3}). On note $\mathcal{T}$ l'adh\'erence sch\'ematique de $T$ dans $\mathcal{G}$. C'est un sch\'ema en groupes de type fini sur $V$. On note $\mathcal{G}^{\textup{ss}}$ le sous-groupe d\'eriv\'e de $\mathcal{G}$, et $\mathcal{G}^{\textup{sc}}$ le rev\^etement simplement connexe de $\mathcal{G}^{\textup{ss}}$. On note aussi $\mathcal{Z}$ le centre de $\mathcal{G}$ et $\mathcal{Z}^{\textup{sc}}$ celui de $\mathcal{G}^{\textup{sc}}$, et $\mathcal{\rho} : \mathcal{G}^{\textup{sc}} \rightarrow \mathcal{G}$. On utilise alors le corollaire 6.3 de l'expos\'e XV de \cite{SGA3} pour savoir que, quitte \`a r\'eduire $V$, on peut supposer que $\mathcal{T} \subset \mathcal{G}$ est un sous-tore maximal de $\mathcal{G}$ (au sens de \cite{SGA3}). On a alors un diagramme commutatif de sch\'emas en groupes sur $V$ (dont les lignes sont les immersions ferm\'ees naturelles) :
\begin{displaymath}
\xymatrix{
\mathcal{Z}^{\textup{sc}} \ar[r] \ar[d]^{\rho} & \mathcal{T}^{\textup{sc}} \ar[r] \ar[d]^{\rho} & \mathcal{G}^{\textup{sc}} \ar[d]^{\rho} \\
\mathcal{Z} \ar[r] & \mathcal{T} \ar[r] & \mathcal{G} \\
}
\end{displaymath}

On dispose ainsi de trois modules crois\'es sur $V$, qui sont
quasi-isomorphes (voir \cite{BorAMS}, lemme 3.8.1) :
$[\mathcal{G}^{\textup{sc}} \rightarrow \mathcal{G}]$,
$[\mathcal{Z}^{\textup{sc}} \rightarrow \mathcal{Z}]$ et $\mathcal{C}
:= [\mathcal{T}^{\textup{sc}} \rightarrow \mathcal{T}]$. Ces modules
crois\'es permettent de d\'efinir un morphisme naturel $\textup{ab}^0
: H^0(V, \mathcal{G}) \rightarrow \H^0(V, \mathcal{C})$ (pour des
pr\'ecisions sur ce morphisme d'ab\'elianisation, voir \cite{Del},
\cite{Bre} ou \cite{BorAMS}).

On note aussi, pour un ensemble fini de places $S$ suffisamment
grand (i.e. contenant $S_1$), $\mathcal{P}^i_{S}(\mathcal{G}) :=
\prod_{v \in S} H^i(k_v, G) \times \prod_{v \notin S}
H^i(\mathcal{O}_v, \mathcal{G})$ pour $i = 0$ ou $1$.

\begin{lem}
$\textup{Ker}(\rho : G^{\textup{sc}} \rightarrow G)$ est un groupe alg\'ebrique ab\'elien fini.
\end{lem}

\begin{dem}
En effet, $\textup{ker}(\rho) = \textup{Ker}( G^{\textup{sc}} \rightarrow G^{\textup{ss}})$ est le groupe fondamental de $G^{\textup{ss}}$ (voir \cite{BorAMS}, lemme 2.4.1)
\end{dem}

On va avoir besoin ici du lemme suivant, qui dit essentiellement que la compl\'etion profinie commute au produit fini de groupes topologiques :

\begin{lem}
\label{lem compl prod}
Soient $G$ et $H$ deux groupes topologiques. Alors il existe un isomorphisme canonique $(G \times H)^{\wedge} \cong G^{\wedge} \times H^{\wedge}$.
\end{lem}

\begin{dem}
On dispose d'un morphisme canonique $\phi : (G \times H)^{\wedge} \rightarrow G^{\wedge} \times H^{\wedge}$.
\begin{itemize}
	\item $\phi$ est injectif :
pour montrer cette injectivit\'e, il suffit de montrer que tout sous-groupe distingu\'e $K$ ouvert d'indice fini dans $G \times H$ contient un produit $K_1 \times K_2$, o\`u $K_1$ (resp. $K_2$) est un sous-groupe distingu\'e ouvert d'indice fini de $G$ (resp. $H$). Soit donc un tel sous-groupe $K$. On note $G'$ (resp. $H'$) le sous-groupe $G \times 1$ (resp. $1 \times H$) de $G \times H$; ce sous-groupe est isomorphe comme groupe topologique \`a $G$ (resp. $H$).
On dispose du diagramme commutatif exact de groupes topologiques suivant :
\begin{displaymath}
\xymatrix{
1 \ar[r] & G' \cap K \ar[r] \ar[d] & G' \ar[r] \ar[d] & G' / (G' \cap K) \ar[d] \ar[r] & 1 \\
1 \ar[r] & K \ar[r] & G \times H \ar[r] & (G \times H) / K \ar[r] & 1
}
\end{displaymath}
o\`u la derni\`ere fl\`eche verticale est injective.
Alors $G' \cap K$ est un sous-groupe distingu\'e ouvert d'indice fini de $G'$, et d\'efinit naturellement un sous-groupe distingu\'e ouvert d'indice fini $K_1$ de $G$, de sorte que $K_1 \times 1 \subset K$. De m\^eme on construit $K_2 \subset H$ sous-groupe ouvert d'indice fini tel que $1 \times K_2 \subset K$. Alors il est clair que $K_1 \times K_2$ est contenu dans $K$, ce qui assure l'injectivit\'e de $\phi$.
	\item $\phi$ est surjectif :
pour la surjectivit\'e, on remarque que par fonctorialit\'e, les morphismes $p_1 : (G \times H)^{\wedge} \rightarrow G^{\wedge}$ et $p_2 : (G \times H)^{\wedge} \rightarrow H^{\wedge}$ (induits par les projections $G \times H \rightarrow G, H$) admettent des sections $i_1$ et $i_2$. Or $\phi(x) = (p_1(x), p_2(x))$ pour tout $x \in (G \times H)^{\wedge}$, donc tout \'el\'ement $(a,b) \in G^{\wedge} \times H^{\wedge}$ s'\'ecrit $(a,b) = (a,1).(1,b) = \phi(i_1(a)).\phi(i_2(b)) = \phi(i_1(a).i_2(b))$, donc $\phi$ est surjectif. 
\end{itemize}
\end{dem}

D\'esormais, $U = \textup{Spec}(\mathcal{O}_{k, S})$ est un ouvert de $V$, avec $S_1 \subset S$.
On va maintenant utiliser les morphismes d'ab\'elianisation pour le $H^0$. Montrons d'abord que le conoyau du morphisme d'ab\'elianisation est control\'e par le $H^1$ du rev\^etement universel de $G$ :

\begin{lem}
\label{coker ab}
Soit $X$ un sch\'ema, $H$ un $X$-sch\'ema en groupes r\'eductif.
Alors $\textup{ab}^0 : H^0(X, H) \rightarrow \H^0_{\textup{ab}}(X, H)$ est un morphisme de groupes, de conoyau s'injectant dans $H^1(X, H^{\textup{sc}})$ de sorte que le diagramme suivant commute ($\delta$ est le morphisme d'ensembles point\'es provenant du morphisme de modules crois\'es $[H^{\textup{sc}} \rightarrow H] \rightarrow [H^{\textup{sc}} \rightarrow 1]$) :
\begin{displaymath}
\xymatrix{
H^0(X, H) \ar[r]^{\textup{ab}^0} & H^0_{\textup{ab}}(X, H) \ar[rd] \ar[rr]^{\delta} & & H^1(X, H^{\textup{sc}}) \\
& & \textup{Coker}(\textup{ab}^0) \ar[ru]^{i} &
}
\end{displaymath}
\end{lem}

\begin{dem}
On sait que $\textup{ab}^0$ est un morphisme de groupes (voir par exemple \cite{Bre}). Il reste \`a montrer que l'application $i$ du diagramme pr\'ec\'edent est bien d\'efinie. Pour cela, il suffit de montrer que pour tout $g' \in H^0_{\textup{ab}}(X, H)$ et tout $g \in H^0(X, H)$, on a $\delta(g'.\textup{ab}^0(g)) = \delta(g')$. Mais ceci est \'evident par d\'efinition de la structure de groupe sur $H^0_{\textup{ab}}(X, H)$ et du morphisme $\delta$ (voir par exemple \cite{BorAMS} 3.3.1, pour une d\'efinition en termes de cocycles).
\end{dem}

\begin{lem}
\label{lem sel compl}
On a une suite exacte naturelle d'ensembles point\'es (dont les trois
premiers termes forment une suite exacte de groupes) :
$$\mathcal{P}^0_S(\mathcal{G}^{\textup{sc}})^{\wedge} \xrightarrow{\rho} \mathcal{P}^0_S(\mathcal{G})^{\wedge} \xrightarrow{\textup{ab}^0} \mathcal{P}^0_S(\mathcal{C})^{\wedge} \rightarrow \mathcal{P}^1_S(\mathcal{G}^{\textup{sc}})$$
\end{lem}

\begin{dem}
On remarque d'abord que l'ensemble point\'e
$\mathcal{P}_S^1(\mathcal{G}^{\textup{sc}})$ est fini et est isomorphe
\`a $\prod_{v \in \Omega_{\infty}} H^1(k_v, G^{\textup{sc}})$ (voir
\cite{PR}, th\'eor\`eme 6.4 pour les places $v \in S \setminus S_{\infty}$, ainsi que le th\'eor\`eme de Lang pour
les places $v \notin S$). 
Montrons que la suite 
$$\mathcal{P}^0_S(\mathcal{G})^{\wedge} \xrightarrow{\textup{ab}^0} \mathcal{P}^0_S(\mathcal{C})^{\wedge} \rightarrow \mathcal{P}^1_S(\mathcal{G}^{\textup{sc}})$$
est exacte ($\mathcal{P}^1_S(\mathcal{G}^{\textup{sc}})$ est un ensemble point\'e fini, les autres ensembles \'etant des groupes).
Si $v$ est une place finie hors de $S$, le morphisme d'ab\'elianisation $\textup{ab}^0_v : H^0(\mathcal{O}_v, \mathcal{G}) \rightarrow \H^0(\mathcal{O}_v, \mathcal{C})$ est un morphisme surjectif de groupes compacts (puisque $H^1(\mathcal{O}_v, \mathcal{G}^{\textup{sc}}) = 1$ par le th\'eor\`eme de Lang), donc la surjection $\prod_{v \textup{ finie } \notin S} H^0(\mathcal{O}_v, \mathcal{G}) \rightarrow \prod_{v \textup{ finie } \notin S} \H^0(\mathcal{O}_v, \mathcal{C})$ n'est pas affect\'ee par l'op\'eration de compl\'etion profinie. Si $v$ est une place finie dans $S$, la compl\'etion de la surjection $H^0(k_v, G) \rightarrow \H^0(k_v, C)$ donne une surjection $H^0(k_v, G)^{\wedge} \rightarrow \H^0(k_v, C)^{\wedge}$. Enfin, si $v$ est une place r\'eelle, on a une suite exacte de groupes finis
$$\widehat{H}^0(\R, G) \rightarrow \widehat{H}^0(\R, C) \rightarrow Q_v \rightarrow 0$$
et une suite exacte d'ensembles point\'es
$$\widehat{H}^0(\R, G) \rightarrow \widehat{H}^0(\R, C) \rightarrow H^1(\R, G^{\textup{sc}})$$
avec $Q_v$ s'injectant naturellement dans $H^1(\R, G^{\textup{sc}})$ par le lemme \ref{coker ab}. Donc en prenant la compl\'etion profinie de la premi\`ere suite, on obtient la m\^eme suite exacte, et on peut bien envoyer $Q_v$ dans $H^1(\R, G^{\textup{sc}})$ de mani\`ere injective. Ainsi, \`a l'aide du lemme \ref{lem compl prod}, on conclut que la suite 
$$\mathcal{P}^0_S(\mathcal{G})^{\wedge} \xrightarrow{\textup{ab}^0} \mathcal{P}^0_S(\mathcal{C})^{\wedge} \rightarrow \mathcal{P}^1_S(\mathcal{G}^{\textup{sc}})$$
est exacte.

Reste \`a montrer l'exactitude de la suite de l'\'enonc\'e en $\mathcal{P}^0_S(\mathcal{G})^{\wedge}$.
Puisque la compl\'etion profinie commute au produit fini (voir lemme \ref{lem compl prod}), il suffit de montrer l'exactitude au niveau de chaque place $v \in S$ : en effet, la suite exacte $\prod_{v \notin S} \mathcal{G}^{\textup{sc}}(\mathcal{O}_v) \rightarrow \prod_{v \notin S} \mathcal{G}(\mathcal{O}_v) \rightarrow \prod_{v \notin S} \H^0(\mathcal{O}_v, \mathcal{C}) \rightarrow 0$ est une suite exacte de groupes compacts ($H^1(\mathcal{O}_v, \mathcal{G}^{\textup{sc}}) = 0$ gr\^ace au th\'eor\`eme de Lang), donc sa compl\'etion profinie est \'egale \`a elle m\^eme. Par cons\'equent, en voyant $\mathcal{P}^0_S(\mathcal{G})^{\wedge}$ comme le produit fini de $\prod_{v \notin S} \mathcal{G}(\mathcal{O}_v)$ et des $H^0(k_v, G)$ pour $v \in S$, le lemme \ref{lem compl prod} assure bien qu'il suffit de montrer l'exactitude de le suite compl\'et\'ee en chaque place $v \in S$. Si $v \in S$ est non archim\'edienne, on consid\'ere la suite exacte courte suivante :
$$1 \rightarrow \textup{Ker}(\textup{ab}^0_v) \rightarrow G(k_v) \xrightarrow{\textup{ab}^0_v} \textup{Im}(\textup{ab}^0_v) \rightarrow 0$$
Or le morphisme $\textup{ab}^0_v : G(k_v) \rightarrow \H^0(k_v,C)$ est ouvert, par d\'efinition de la topologie sur $\H^0(k_v, C)$ et par commutativit\'e du diagramme 
\xymatrix{
G(k_v) \ar[r]^{\textup{ab}^0_v} & \H^0(k_v, C) \\
T(k_v) \ar[u] \ar[ru] &
}
Donc la compl\'etion profinie fournit une suite exacte
$$\textup{Ker}(\textup{ab}^0_v)^{\wedge} \rightarrow G(k_v)^{\wedge} \xrightarrow{\textup{ab}^0} \textup{Im}(\textup{ab}^0_v)^{\wedge}$$
D'o\`u un diagramme commutatif :
\begin{displaymath}
\xymatrix{
G^{\textup{sc}}(k_v)^{\wedge} \ar[r] \ar[d] & G(k_v)^{\wedge} \ar[r] \ar[d]^{=} & \H^0(k_v, C)^{\wedge} \\
\textup{Ker}(\textup{ab}^0_v)^{\wedge} \ar[r] & G(k_v)^{\wedge} \ar[r] & \textup{Im}(\textup{ab}^0_v)^{\wedge} \ar[u]
}
\end{displaymath}
Or la premi\`ere fl\`eche verticale est surjective, et la derni\`ere est injective car $\textup{Im}(\textup{ab}^0_v)$ est un sous-groupe de $\H^0(k_v, C)$ d'indice fini. Une chasse au diagramme assure alors l'exactitude de la premi\`ere ligne de ce diagramme.
Pour les places archim\'ediennes dans $S$, le raisonnement est le m\^eme en utilisant les groupes de composantes connexes.
\end{dem}

\paragraph{Obstruction de Brauer-Manin et ab\'elianisation}
On commence par rappeler la d\'efinition de l'\emph{accouplement de Brauer-Manin} sur $G$ (voir par exemple le livre de Skorobogatov \cite{Sko}, section 5.2 pour davantage de d\'etails). On consid\`ere l'application suivante :
\fct{P^0(k, G) \times \textup{Br}(G)}{\Q / \Z}{(P, A)}{\langle A, P \rangle_{\textup{BM}} := \sum_{v \in \Omega_k} j_v(A(P_v))}
o\`u $j_v : \textup{Br}(k_v) \rightarrow \Q / \Z$ est l'invariant donn\'e par la th\'eorie du corps de classes local et $A(P_v) \in \textup{Br}(k_v)$ est l'\'evaluation de $A \in \textup{Br}(G)$ en $P_v \in H^0(k_v, G)$. On peut alors d\'efinir, pour toute partie $B$ de $\textup{Br}(G)$ l'ensemble suivant :
$$P^0(k, G)^{B} := \left\{ P \in P^0(k, G) : \langle A, P \rangle_{\textup{BM}} = 0, \forall A \in B \right\}$$
Gr\^ace \`a la loi de r\'eciprocit\'e de la th\'eorie du corps de classes global, on sait que $\overline{G(k)} \subset P^0(k, G)^{B}$. Enfin, les \'el\'ements de $\textup{Br } k \subset \textup{Br } G$ sont orthogonaux \`a tous les points de $P^0(k, G)$, donc l'accouplement pr\'ec\'edent se factorise en une application $P^0(k, G) \times \textup{Br}(G) / \textup{Br } k \rightarrow \Q / \Z$, et fournit donc un morphisme $P^0(k, G) \rightarrow ( \textup{Br}(G) / \textup{Br } k)^D$. 

Or, gr\^ace aux travaux de Borovoi et Van Hamel, on sait que l'on dispose d'un isomorphisme canonique $\kappa : \H^1(k, \widehat{C}) \cong \textup{Br}_a(G)$ (voir \cite{BVH}, th\'eor\`eme 4.8 et corollaire 2.20.(ii)). On note \'egalement $\kappa$ l'isomorphisme correspondant par dualit\'e $\kappa : \textup{Br}_a(G)^D \cong \H^1(k, \widehat{C})^D$.

On consid\`ere d\'esormais le diagramme suivant :
\begin{eqnarray}
\label{diag abel}
\xymatrix{
H^0(U, \mathcal{G}^{\textup{sc}})^{\wedge} \ar[r] \ar[d] & \mathcal{P}^0_S(\mathcal{G}^{\textup{sc}})^{\wedge} \ar[d] & \\
H^0(U, \mathcal{G})^{\wedge} \ar[r] \ar[d] & \mathcal{P}_S^0(\mathcal{G})^{\wedge} \ar[r]^{\theta_{BM}} \ar[d] & \textup{Br}_a(G)^D \ar[d]^{\kappa} \\
\H^0(U, \mathcal{C})^{\wedge} \ar[r] \ar[d] & \mathcal{P}_S^0(\mathcal{C})^{\wedge} \ar[r]^{\theta} \ar[d] & \H^1(k, \widehat{C})^D \\
H^1(U, \mathcal{G}^{\textup{sc}}) \ar[r] & \mathcal{P}^1_S(\mathcal{G}^{\textup{sc}}) &
}
\end{eqnarray}
les fl\`eches sans nom dans ce diagramme \'etant les fl\`eches \'evidentes, \`a l'exception de la fl\`eche $\H^0(U, \mathcal{C})^{\wedge} \rightarrow H^1(U, \mathcal{G}^{\textup{sc}})$, qui est obtenue \`a partir de la fl\`eche $\H^0(U, \mathcal{C}) \rightarrow H^1(U, \mathcal{G}^{\textup{sc}})$ gr\^ace au lemme \ref{coker ab} et \`a la finitude de $H^1(U, \mathcal{G}^{\textup{sc}})$. La fl\`eche $\theta_{BM}$ est la fl\`eche induite par l'accouplement de Brauer-Manin, et la fl\`eche $\theta$ est la fl\`eche induite par la dualit\'e locale pour les complexes de tores.

On va d'abord montrer que le carr\'e de droite dans ce diagramme est commutatif, \`a savoir le lemme suivant :

\begin{lem}
\label{BMlocal}
Soit $\theta_{BM} : P^0(k, G) \rightarrow \textup{Br}_a(G)^D$ induit par l'accouplement de Brauer-Manin, et $\theta : \P^0(k, C) \rightarrow \H^1(k, \widehat{C})^D$ l'accouplement induit par la dualit\'e locale. Alors le diagramme 
\begin{displaymath}
\xymatrix{
P^0(k, G) \ar[d]^{\textup{ab}^0} \ar[r]^{\theta_{BM}} & \textup{Br}_a(G)^D \ar[d]^{\kappa} \\
\P^0(k, C) \ar[r]^{\theta} & \H^1(k, \widehat{C})^D
}
\end{displaymath}
est commutatif (au signe pr\`es) et fonctoriel en $C$.
\end{lem}

\begin{dem}
Les deux fl\`eches $\theta_{BM}$ et $\theta$ \'etant d\'efinies par des sommes de contributions locales, il suffit de montrer la commutativit\'e du diagramme suivant :
\begin{displaymath}
\xymatrix{
H^0(k_v, G) \ar[d]^{\textup{ab}_v^0} \ar[r]^{\theta_{BM, v}} & \textup{Br}_1(G_v)^D \ar[d]^{\kappa'_v} \\
\H^0(k_v, C) \ar[r]^{\theta_v} & \H^1(k_v, \widehat{C})^D
}
\end{displaymath}
o\`u $\kappa'_v$ est obtenu en dualisant la compos\'ee de $\kappa_v : \H^1(k_v, \widehat{C}) \rightarrow \textup{Br}_a(G_v)$ avec l'inclusion $\textup{Br}_a(G_v) \rightarrow \textup{Br}_1(G_v)$ donn\'ee par la section unit\'e de $G_v$.
Montrons qu'il suffit alors de consid\'erer seulement le cas o\`u $G$ est un $k$-tore : on sait qu'il existe un $k$-groupe r\'eductif $H$, avec $H^{\textup{ss}}$ simplement connexe, et un morphisme surjectif $\phi : H \rightarrow G$, dont le noyau est un tore quasi-trivial (voir par exemple \cite{MS}, p 297, proposition 3.1). Supposons que pour $H$, le diagramme local soit commutatif. Par fonctorialit\'e de la fl\`eche d'abelianisation $\textup{ab}^0$ (voir \cite{BorAMS}, proposition 3.11), un tel $H$ induit le diagramme suivant (o\`u $C'$ est le complexe de tores associ\'e \`a $H$) :
\begin{displaymath}
\xymatrix{
H(k_v) \ar[dd]^{\phi} \ar[rd]^{\textup{ab}^0_v} \ar[rr]^{\theta_{BM, v}}  & & \textup{Br}_1(H_v)^D \ar[rd]^{\kappa'_v} \ar'[d][dd]^(.25){\phi_*} & \\
& \H^0(k_v, C') \ar[dd]^(.3){\phi'} \ar[rr]^(.4){\theta_v} & & \H^1(k_v, \widehat{C'})^D \ar[dd]^{\phi'_*} \\
G(k_v) \ar[rd]_{\textup{ab}^0_v} \ar'[r]^(.7){\theta_{BM, v}}[rr] & & \textup{Br}_1(G_v)^D \ar[rd]^{\kappa'_v} & \\
& \H^0(k_v, C) \ar[rr]^{\theta_v} & & \H^1(k_v, \widehat{C})^D
}
\end{displaymath}
o\`u toutes les faces verticales du cube sont des diagrammes commutatifs (fonctorialit\'es de l'ab\'elianisation, de l'accouplement de Brauer-Manin, de la dualit\'e locale et de l'isomorphisme $\kappa_v$). Or la fl\`eche $\phi : H(k_v) \rightarrow G(k_v)$ est surjective, donc la commutativit\'e de la face sup\'erieure et des faces verticales impliquent la commutativit\'e de la face inf\'erieure. Par cons\'equent, il suffit de prouver le r\'esultat pour $G$ groupe r\'eductif tel que $G^{\textup{ss}}$ est simplement connexe. On consid\`ere dans ce cas la suite exacte suivante :
$$1 \rightarrow G^{\textup{ss}} \rightarrow G \rightarrow T \rightarrow 1$$
o\`u $T$ est un $k$-tore. Supposons que l'on sache d\'emontrer la commutativit\'e du diagramme pour le tore $T$. D\'eduisons-en la propri\'et\'e pour le groupe $G$ : on dispose du diagramme suivant, dont les colonnes sont exactes :

\begin{displaymath}
\xymatrix{
G^{\textup{ss}}(k_v) \ar[dd] \ar[rd]^{\textup{ab}^0_v} \ar[rr]^{\theta_{BM, v}}  & & \textup{Br}_1(G^{\textup{ss}}_v)^D \ar[rd]^{\kappa'_v} \ar'[d][dd] & \\
& \H^0(k_v, C') \ar[dd] \ar[rr]^(.4){\theta_v} & & \H^1(k_v, \widehat{C'})^D \ar[dd] \\
G(k_v) \ar[dd] \ar[rd]_{\textup{ab}^0_v} \ar'[r]^(.7){\theta_{BM, v}}[rr] & & \textup{Br}_1(G_v)^D \ar'[d][dd] \ar[rd]^{\kappa'_v} & \\
& \H^0(k_v, C) \ar[dd] \ar[rr]^(0.4){\theta_v} & & \H^1(k_v, \widehat{C})^D \ar[dd] \\
T(k_v)  \ar[rd]_{\textup{ab}^0_v} \ar'[r]^(.7){\theta_{BM, v}}[rr] & & \textup{Br}_1(T_v)^D \ar[rd]^{\kappa'_v} & \\
& H^0(k_v, T) \ar[rr]^{\theta_v} & & H^2(k_v, \widehat{T})^D
}
\end{displaymath}
Or par fonctorialit\'e les faces verticales sont des diagrammes commutatifs; la face inf\'erieure est un diagramme commutatif par hypoth\`ese. Enfin, la nullit\'e de $\textup{Br}_a(G^{\textup{ss}}_v) \cong \H^1(k_v, \widehat{C'})$ assure que le morphisme $\textup{Br}_a(G_v)^D = \H^1(k_v, \widehat{C})^D \rightarrow \textup{Br}_a(T_v)^D = H^2(k_v, \widehat{T})$ est injectif (voir par exemple \cite{San}, corollaire 6.11.2). Une chasse au diagramme imm\'ediate assure alors la commutativit\'e du carr\'e horizontal central. Par cons\'equent, on s'est ramen\'e \`a montrer la propri\'et\'e de compatibilit\'e au cas o\`u $G$ est un tore. Or ce cas est connu (voir par exemple \cite{HSz2}, preuve du th\'eor\`eme 6.1). Par cons\'equent, le lemme \ref{BMlocal} est d\'emontr\'e.
\end{dem}

Par cons\'equent, le diagramme (\ref{diag abel}) est commutatif (au signe pr\`es).

On conna\^it en outre les faits suivants : 
la deuxi\`eme colonne de ce diagramme est exacte (c'est le lemme \ref{lem sel compl}), la premi\`ere est exacte en $\H^0(U, \mathcal{C})^{\wedge}$ et la troisi\`eme ligne est exacte (voir proposition \ref{prop ouvert compl}).

\begin{lem}
\label{Nisnevich}
On suppose que $G^{\textup{sc}}$ v\'erifie l'approximation hors de $S_0$, pour $S_0 \subset S$. Alors l'application $H^1(U, \mathcal{G}^{\textup{sc}}) \rightarrow  \mathcal{P}^1_S(\mathcal{G}^{\textup{sc}})$ a un noyau trivial.
\end{lem}

\begin{dem}
Par un r\'esultat de Nisnevich (voir \cite{Nis}, th\'eor\`eme 2.1)
d\'emontr\'e par Gille (th\'eor\`eme 5.1 dans l'appendice de \cite{Gil}), si on note $c(G) := G(\A_S) \backslash G(\A_k) / G(k)$ le groupe des $S$-classes de $G$, on dispose d'une suite exacte
$$1 \rightarrow c(G^{\textup{sc}}) \rightarrow H^1(U, \mathcal{G}^{\textup{sc}}) \rightarrow H^1(k, G^{\textup{sc}})$$
puisque par th\'eor\`eme de Lang, les ensembles $H^1(\mathcal{O}_v,
G^{\textup{sc}})$ sont triviaux pour $v \notin S$. Or le groupe
$G^{\textup{sc}}$ v\'erifie l'approximation forte hors de $S_0 \subset S$, donc l'ensemble $c(G)$ est trivial (voir \cite{PR}, proposition 5.4). Cela assure le r\'esultat, puisque le morphisme $H^1(k, G^{\textup{sc}}) \rightarrow \P^1(k, G^{\textup{sc}})$ est injectif (principe de Hasse pour les groupes semi-simples simplement connexes, voir par exemple \cite{PR}, th\'eor\`eme 6.6).
\end{dem}

Introduisons d\'esormais le d\'efaut d'approximation forte sur $U$, \`a savoir le quotient topologique :
$$A_{S, S_0}(G) := \mathcal{P}^0_{S}(\mathcal{G}) / \overline{H^0(U,
  \mathcal{G}).G^{\textup{sc}}_{S_0}}$$
On dispose alors du r\'esultat suivant (le cas semi-simple est la
proposition 8.8 de \cite{PR}) :

\begin{prop}
Soit $G$ un groupe r\'eductif sur $k$.On suppose que $G^{\textup{sc}}$
v\'erifie l'approximation forte hors de $S_0 \subset S$.
Quitte \`a r\'eduire $U$, l'ensemble $\overline{H^0(U, \mathcal{G}).G^{\textup{sc}}_{S_0}}$ est un sous-groupe distingu\'e de $\mathcal{P}^0_{S}(\mathcal{G})$ et le quotient $A_{S, S_0}(G)$ est un groupe topologique ab\'elien.
\end{prop}

\begin{dem}
Quitte \`a r\'eduire $U$, on peut supposer $\mathcal{G}$ r\'eductif sur $U$.
On consid\`e le diagramme commutatif suivant :
\begin{displaymath}
\xymatrix{
\overline{H^0(U, \mathcal{G}^{\textup{sc}}).G^{\textup{sc}}_{S_0}} \ar[r] \ar[d]^{\rho} & \mathcal{P}^0_S(\mathcal{G}^{\textup{sc}}) \ar[d]^{\rho} \\
\overline{H^0(U, \mathcal{G}).\rho(G^{\textup{sc}}_{S_0})} \ar[r] \ar[d] & \mathcal{P}^0_S(\mathcal{G}) \ar[d]^{\textup{ab}^0} \\
\overline{\H^0(U, \mathcal{C})} \ar[r] & \mathcal{P}^0_S(\mathcal{C})
}
\end{displaymath}
On sait que $\overline{\H^0(U, \mathcal{C})}$ est un sous-groupe du groupe ab\'elien $\mathcal{P}^0_S(\mathcal{C})$, et que la premi\`ere fl\`eche horizontale est une \'egalit\'e (par le th\'eor\`eme d'approximation forte). Soit alors $h \in \overline{H^0(U, \mathcal{G}).\rho(G^{\textup{sc}}_{S_0})}$ et $g \in \mathcal{P}^0_S(\mathcal{G})$. Montrons que l'\'el\'ement $h' := g h g^{-1} h^{-1}$ est dans $\overline{H^0(U, \mathcal{G}).\rho(G^{\textup{sc}}_{S_0})}$. On a $\textup{ab}^{0}(h') = 0$ car $\mathcal{P}^0_S(\mathcal{C})$ est ab\'elien. Donc $h'$ se rel\`eve en un \'el\'ement $h''$ dans $\mathcal{P}^0_S(\mathcal{G}^{\textup{sc}})$. Par th\'eor\`eme d'approximation forte, $h'' \in \overline{H^0(U, \mathcal{G}^{\textup{sc}}).G^{\textup{sc}}_{S_0}}$, et donc $h' = \rho(h'') \in \overline{H^0(U, \mathcal{G}).\rho(G^{\textup{sc}}_{S_0})}$, ce qui assure que $\overline{H^0(U, \mathcal{G}).\rho(G^{\textup{sc}}_{S_0})}$ est un sous-groupe distingu\'e de $\mathcal{P}^0_S(\mathcal{G})$. Montrons d\'esormais que le quotient $Q := \mathcal{P}^0_S(\mathcal{G}) / \overline{H^0(U, \mathcal{G}).\rho(G^{\textup{sc}}_{S_0})}$ est ab\'elien : soient $a, b \in Q$, notons $a'$, $b'$ des relev\'es de $a$ et $b$ dans $\mathcal{P}^0_S(\mathcal{G})$. Alors $g := a' b' a'^{-1} b'^{-1}$ s'envoie sur $0$ dans le groupe ab\'elien $\mathcal{P}^0_S(\mathcal{C})$, donc comme plus haut, $g$ se rel\`eve en $h \in \mathcal{P}^0_S(\mathcal{G}^{\textup{sc}})$. Or par approximation forte, $h \in \overline{H^0(U, \mathcal{G}^{\textup{sc}}).G^{\textup{sc}}_{S_0}}$, donc $g = \rho(h) \in \overline{H^0(U, \mathcal{G}).\rho(G^{\textup{sc}}_{S_0})}$, ce qui assure que le commutateur $a b a^{-1} b^{-1}$ est trivial dans $Q$, donc $Q$ est ab\'elien. 
\end{dem}
Montrons alors le r\'esultat suivant :
\begin{theo} 
\label{theo compl red}
Soit $S_0 \subset S$ un ensemble fini de places de $k$ tel que $G^\textup{sc}$ v\'erifie la propri\'et\'e d'approximation forte hors de $S_0$. Alors la suite suivante
$$\left( \overline{H^0(U, \mathcal{G}) . \rho(G^{\textup{sc}}_{S_0})}
\right)^{\wedge} \rightarrow \mathcal{P}^0_S(\mathcal{G})^{\wedge}
\rightarrow (\textup{Br}_a G)^D$$ est exacte et fonctorielle en G.
\end{theo}

\begin{dem}
Cette suite est un complexe car les \'el\'ements de $G^{\textup{sc}}_{S_0}$ ont une image nulle dans $(\textup{Br}_a G)^D$ (voir le diagramme commutatif (\ref{diag abel})). \\
Soit $(g_v) \in \mathcal{P}^0_S(\mathcal{G})^{\wedge}$ d'image nulle dans $(\textup{Br}_a G)^D$. On note $(g'_v)$ son image dans $\mathcal{P}^0_S(\mathcal{C})^{\wedge}$. Par exactitude de la troisi\`eme ligne du diagramme (\ref{diag abel}), $(g'_v)$ se rel\`eve en un \`el\'ement $g' \in \H^0(U, \mathcal{C})^{\wedge}$ qui s'envoie sur $1$ dans $H^1(U, \mathcal{G}^{\textup{sc}})$. Par exactitude de la premi\`ere colonne, $g'$ se rel\`eve donc en un \'el\'ement $\tilde{g} \in H^0(U, \mathcal{G})^{\wedge}$ dont l'image $(\tilde{g}_v)$ dans $\mathcal{P}^0_S(\mathcal{G})^{\wedge}$ a m\^eme image que $(g_v)$ dans $\mathcal{P}^0_S(\mathcal{C})^{\wedge}$. Donc par exactitude de la deuxi\`eme colonne, l'\'el\'ement $(\tilde{g}_v)^{-1}.(g_v) \in \mathcal{P}^0_S(\mathcal{G})^{\wedge}$ se rel\`eve en un \`el\'ement $(h_v) \in \mathcal{P}^0_S(\mathcal{G}^{\textup{sc}})^{\wedge}$.
On utilise d\'esormais le lemme suivant :

\begin{lem}
La fl\`eche naturelle $\left( \overline{H^0(U, \mathcal{G}^{\textup{sc}}). G^{\textup{sc}}_{S_0}} \right)^{\wedge} \rightarrow \mathcal{P}^0_S(\mathcal{G}^{\textup{sc}})^{\wedge}$ est surjective.
\end{lem}

\begin{dem}
C'est une cons\'equence imm\'ediate de la propri\'et\'e d'approximation forte pour le groupe semi-simple simplement connexe $G^{\textup{sc}}$ hors des places $S_0$.
\end{dem}

Avec ce lemme, l'\'el\'ement $(h_v)$ provient d'un \'el\'ement $h \in \left( \overline{H^0(U, \mathcal{G}^{\textup{sc}}). G^{\textup{sc}}_{S_0}} \right)^{\wedge}$. On consid\`ere alors l'image $\rho(h) \in \left( \overline{H^0(U, \mathcal{G}). \prod_{v \in S_0} H^0(k_{v}, G) } \right)^{\wedge}$ : par fonctorialit\'e, l'\'el\'ement $\rho(h).\tilde{g} \in \left(\overline{H^0(U, \mathcal{G}).G^{\textup{sc}}_{S_0}} \right)^{\wedge}$ s'envoie alors sur $(g_v)$ dans $\mathcal{P}^0_S(\mathcal{G})^{\wedge}$.
\end{dem}

\begin{cor}
\label{coro AF ouvert}
On a une suite exacte fonctorielle en $G$:
$$1 \rightarrow \overline{H^0(U, \mathcal{G}).G^{\textup{sc}}_{S_0}} \rightarrow \mathcal{P}^0_S(\mathcal{G}) \rightarrow (\textup{Br}_a G)^D$$
\end{cor}

\begin{dem}
On consid\`ere le conoyau topologique $Q = A_{S, S_0}(G)$ de l'injection $\overline{H^0(U, \mathcal{G}).G^{\textup{sc}}_{S_0}} \rightarrow \mathcal{P}^0_S(\mathcal{G})$. On dispose alors du diagramme commutatif suivant :
\begin{displaymath}
\xymatrix{
0 \ar[r] & \overline{H^0(U, \mathcal{G}).G^{\textup{sc}}_{S_0}} \ar[r] \ar[d] & \mathcal{P}^0_S(\mathcal{G}) \ar[r] \ar[d] & Q \ar[d] \\
 & \overline{H^0(U, \mathcal{G}).G^{\textup{sc}}_{S_0}}^{\wedge} \ar[r] & \mathcal{P}^0_S(\mathcal{G})^{\wedge} \ar[r] & Q^{\wedge}
}
\end{displaymath}
La premi\`ere ligne est exacte par d\'efinition, la seconde ligne est un complexe, et la derni\`ere fl\`eche verticale est injective car le groupe topologique ab\'elien $Q$, comme quotient de $\mathcal{P}^0_S(\mathcal{G})$, est compactement engendr\'e ($G(k_v)$ est compactement engendr\'e pour $v \in S$). On conclut alors gr\^ace au th\'eor\`eme \ref{theo compl red}, avec une chasse au diagramme (voir \cite{HarAF}, preuve du th\'eor\`eme 1).
\end{dem}

\begin{theo} 
\label{theo fin}
Soit $G$ un $k$-groupe r\'eductif, $S_0$ un ensemble fini de places de $k$ telle que $G^{\textup{sc}}$ v\'erifie l'approximation forte hors de l'ensemble de places $S_0$ (c'est le cas lorsque $(G^{\textup{sc}})_{S_0}^i$ est non compact pour tout $k$-facteur presque $k$-simple $(G^{\textup{sc}})^i$ de $G^{\textup{sc}}$ ). L'ensemble $\textup{\cyr{SH}}^1(k, G)$ est muni d'une structure de groupe ab\'elien via un isomorphisme fonctoriel $\textup{\cyr{B}}(G)^D \cong \textup{\cyr{SH}}^1(k, G)$. Alors l'adh\'erence $\overline{G(k).\rho(G^{\textup{sc}}_{S_0})}$ est un sous-groupe distingu\'e de $P^0(k, G)$ et on a une suite exacte de groupes, fonctorielle en $G$ :
$$1 \rightarrow \overline{G(k).\rho(G^{\textup{sc}}_{S_0})} \rightarrow P^0(k,G) \rightarrow (\textup{Br}_a G)^D \rightarrow \textup{\cyr{SH}}^1(k, G) \rightarrow 0$$
\end{theo}

\begin{dem}
Pour les trois premiers termes, il s'agit de passer \`a la limite sur $U$ dans le corollaire \ref{coro AF ouvert}. Concernant l'exactitude au terme $\textup{Br}_a(G)^D$, il suffit de montrer que les images de $P^0(k, G)$ et de $P^0(k, C)$ dans $\H^1(k, \widehat{C})^D$ co\"incident. Pour cela, on utilise le diagramme commutatif suivant
\begin{displaymath}
\xymatrix{
& G(k) \ar[r] \ar[d] & P^0(k, G) \ar[r] \ar[d] & \textup{Br}_a(G)^D \ar[r] \ar[d]^{\simeq} & \textup{\cyr{SH}}^1(G) \ar[d]^{\simeq} & \\
& \H^0(k, C) \ar[r] \ar[d] & \P^0(k, C) \ar[r] \ar[d] & \H^1(k, \widehat{C})^D \ar[r] & \textup{\cyr{SH}}^1(C) \ar[r] & 0 \\
& H^1(k, G^{\textup{sc}}) \ar[r]^{\simeq} \ar[d] & P^1(k, G^{\textup{sc}}) \ar[d] & & & \\
\textup{\cyr{SH}}^1(G) \ar[r] \ar[d]^{\simeq} & H^1(k, G) \ar[r] \ar[d] & P^1(k, G) \ar[d] & & & \\
\textup{\cyr{SH}}^1(C) \ar[r] & \H^1(k, C) \ar[r] & \P^1(k, C) & & & 	
}
\end{displaymath}
le morphisme d'ab\'elianisation $\textup{\cyr{SH}}^1(G) \rightarrow \textup{\cyr{SH}}^1(C)$ \'etant un isomorphisme par \cite{BorAMS}, th\'eor\`eme 5.12. Soit alors un \'el\'ement $\gamma \in \H^1(k, \widehat{C})^D$ provenant d'un $\alpha \in \P^0(k, C)$. On cherche \`a montrer que ce $\gamma$ se rel\`eve dans $P^0(k, G)$. Soit $\beta \in P^1(k, G^{\textup{sc}})$ l'image de $\alpha$. On rel\`eve alors $\beta$ en un \'el\'ement $\beta \in H^1(k, G^{\textup{sc}})$ par le th\'eor\`eme 6.6 de \cite{PR}. Par commutativit\'e du diagramme et exactitude de la troisi\`eme colonne, l'image $\delta$ de $\beta$ dans $H^1(k, G)$ s'envoie sur $0$ dans $P^1(k, G)$. Donc $\delta \in \textup{\cyr{SH}}^1(G)$. Par exactitude de la deuxi\`eme colonne, $\delta$ s'envoie sur $0$ dans $\H^1(k, C)$, i.e. $\delta$ est dans le noyau de $\textup{\cyr{SH}}^1(G) \rightarrow \textup{\cyr{SH}}^1(C)$. Donc par le th\'eor\`eme 5.12 de \cite{BorAMS}, $\delta$ est trivial, donc par exactitude de la deuxi\`eme colonne, $\beta$ se rel\`eve en $\epsilon \in \H^0(k, C)$. L'image de $\epsilon$ dans $\P^0(k, C)$ d\'efinit un $\mu \in P^0(k, G)$ par exactitude de la troisi\`eme colonne. Alors par commutativit\'e du diagramme, cet \'el\'ement $\mu \in P^0(k, G)$ a m\^eme image que $\alpha$ dans $\textup{Br}_a(G)^D$, \`a savoir $\gamma$, ce qui conclut la preuve.
\end{dem}

On peut sans difficult\'es \'etendre ce r\'esultat au cas des groupes
lin\'eaires connexes quelconques (non n\'ecessairement r\'eductifs) :

\begin{cor}
\label{cor fin}
Soit $G$ un $k$-groupe lin\'eaire connexe, $S_0$ comme plus haut
(i.e. tel que le rev\^etement simplement connexe $G^{\textup{sc}}$ du
quotient r\'eductif $G^{\textup{red}}$ de $G$ v\'erifie
l'approximation forte hors de $S_0$). On note $G^{\textup{scu}}$ le
produit fibr\'e $G \times_{G^{\textup{red}}} G^{\textup{sc}}$. Alors on a une suite exacte de groupes, fonctorielle en $G$ :
$$1 \rightarrow \overline{G(k).\rho(G^{\textup{scu}}_{S_0})} \rightarrow P^0(k,G) \rightarrow (\textup{Br}_a G)^D \rightarrow \textup{\cyr{SH}}^1(k, G) \rightarrow 0$$
\end{cor}

\begin{rem}
{\rm Si le groupe $G$ est commutatif, on peut prendre $S_0 =
  \emptyset$. Si $G$ n'est pas commutatif, $S_0$ doit \^etre non vide.}
\end{rem}

\begin{dem}
On consid\'ere le d\'evissage suivant : 
$$1 \rightarrow R_u G \rightarrow G \rightarrow H \rightarrow 1$$
o\`u $R_u G$ est le radical unipotent de $G$, et $H = G^{\textup{red}}$
un groupe r\'eductif, la suite exacte \'etant scind\'ee. On conclut
\`a partir du th\'eor\`eme \ref{theo fin}  en utilisant la
trivialit\'e de $H^1(k, R_u G)$ et $P^1(k, R_u G)$, ainsi que le fait
bien connu que les groupes unipotents v\'erifient la propri\'et\'e
d'approximation forte, c'est-\`a-dire que pour un $k$-groupe unipotent
$U$, l'adh\'erence de $U(k)$ dans $P^0(k, U)$ est $P^0(k, U)$ tout
entier (en prenant \`a nouveau les groupes de composantes connexes au
niveau des places infinies, et en utilisant le fait que si $v$ est une
place infinie, $U(k_v)$ est connexe).
\end{dem}

\begin{rem}
\label{rem Sansuc}
{\rm
Ce r\'esultat est \`a rapprocher du r\'esultat de Sansuc (voir th\'eor\`eme 8.12 de \cite{San}) qui traite de l'approximation faible dans les groupes connexes, et que l'on peut r\'esumer ainsi : si $G/k$ un groupe lin\'eaire connexe, alors on a une suite exacte de groupes
$$1 \rightarrow \overline{G(k)}^f \rightarrow \prod_v G(k_v) \rightarrow \textup{\cyr{B}}_{\omega}(G)^D \rightarrow \textup{\cyr{SH}}^1(k, G) \rightarrow 0$$
o\`u $ \overline{G(k)}^f $ d\'esigne l'adh\'erence de $G(k)$ dans $ \prod_v G(k_v)$ muni de la topologie produit.
}
\end{rem}

On peut \'egalement reformuler le corollaire \ref{cor fin} \`a la mani\`ere du th\'eor\`eme 8.12 de \cite{San} : pour cela, on introduit le d\'efaut d'approximation forte $A_{S_0}(G) := P^0(k, G) / \overline{G(k).\rho(G^{\textup{scu}}_{S_0})}$, qui est un groupe ab\'elien gr\^ace aux r\'esultats pr\'ec\'edents. On dispose alors de la reformulation suivante :

\begin{cor}
Sous les hypoth\`eses pr\'ec\'edentes, l'accouplement naturel
$$A_{S_0}(G) \times \textup{Br}_a(G) / \textup{\cyr{B}}(G) \rightarrow \Q / \Z$$
est une dualit\'e parfaite de groupes ab\'eliens, fonctorielle en $G$.
\end{cor}

\subsection{Une variante : cas o\`u le groupe v\'erifie l'approximation faible}

Dans cette section, on se propose de montrer l'analogue du
th\'eor\`eme 3 de \cite{HarAF}, ce qui va constituer une sorte de
raffinement du th\'eor\`eme 8.12 de \cite{San}. On rappelle d'abord un
r\'esultat de Sansuc (voir th\'eor\`eme 8.12 de \cite{San} et remarque
\ref{rem Sansuc}) : le groupe $G$ v\'erifie l'approximation faible
\SSI $\textup{\cyr{B}}_{\omega}(G) / \textup{\cyr{B}}(G) = 0$.

\begin{theo}
Soit $G$ un $k$-groupe r\'eductif, $S_0$ un ensemble fini de places de $k$, contenant les places archim\'ediennes, tel que $G$ s'\'etende en un $\textup{Spec}(\mathcal{O}_{k, S_0})$-sch\'ema en groupes r\'eductif $\mathcal{G}$. On suppose que :
\begin{itemize}
	\item $G^{\textup{sc}}$ v\'erifie l'approximation forte hors de $S_0$.
	\item  $G$ v\'erifie l'approximation faible sur $k$, i.e. $\textup{\cyr{B}}_{\omega}(G) / \textup{\cyr{B}}(G) = 0$.
\end{itemize}
Soit $S'$ un ensemble fini de places de $k$ disjoint de $S_0$. Alors il existe $S$ contenant $S_0$ tel que $S \cap S' = \emptyset$ et $\mathcal{G}(\mathcal{O}_{k, S})$ est dense dans $\prod_{v \in S'} \mathcal{G}(\mathcal{O}_v)$.
\end{theo}

\begin{dem}
On commence par le lemme suivant :

\begin{lem}
\label{lem approx faible}
Si $\H^1_{S, S'}(k, \widehat{C})$ d\'esigne le sous-groupe de $\H^1(k, \widehat{C})$ form\'e des \'el\'ements $\alpha$ tels que $\alpha_v$ est orthogonal \`a $G(k_v)$ si $v \in S$ et $\alpha_v$ est orthogonal \`a $\mathcal{G}(\mathcal{O}_v)$ si $v \notin S \cup S'$, on a une suite exacte
$$1 \rightarrow \overline{H^0(U, \mathcal{G})}^{S'} \rightarrow \prod_{v \in S'} \mathcal{G}(\mathcal{O}_v) \rightarrow \left( \H^1_{S, S'}(k, \widehat{C}) / \textup{\cyr{SH}}^1(k, \widehat{C}) \right)^D$$
\end{lem}

\begin{dem}
On consid\`ere la suite exacte suivante (voir th\'eor\`eme \ref{theo fin}) :
$$1 \rightarrow \overline{H^0(k, G).G^{\textup{sc}}_{S_0}} \rightarrow P^0(k, G) \rightarrow \textup{Br}_a(G)^D \rightarrow \textup{\cyr{SH}}^1(G) \rightarrow 0$$
En remarquant que l'adh\'erence de $H^0(U, \mathcal{G}).G^{\textup{sc}}_{S_0}$ dans $\mathcal{P}^0_S(\mathcal{G})$ est exactement l'intersection de $\overline{H^0(k, G).G^{\textup{sc}}_{S_0}}$ avec $\mathcal{P}^0_S(\mathcal{G})$ (dans $P^0(k, G)$), on en d\'eduit l'exactitude de la suite suivante :
$$1 \rightarrow \overline{H^0(U, \mathcal{G}).G^{\textup{sc}}_{S_0}}
\rightarrow \mathcal{P}_{S_0}(\mathcal{G}) \rightarrow
\textup{Br}_a(G)^D$$
On se donne alors $(g_v)_{v \in S'} \in \prod_{v \in S'}
\mathcal{G}(\mathcal{O}_v)$, orthogonal au groupe $\H^1_{S, S'}(k, \widehat{C})$. On plonge $(g_v)$ dans
$\mathcal{P}_{S_0}(\mathcal{G})$ en compl\'etant par l'unit\'e de
$G(k_v)$ pour $v \notin S'$. Alors ce point d\'efinit une obstruction
$\delta \in \textup{Br}_a(G)^D$ qui s'envoie sur $0$ dans $\H^1_{S,
  S'}(k, \widehat{C})^D$. Or on a une suite exacte 
$$\H^1_{S, S'}(k, \widehat{C}) \rightarrow \textup{Br}_a(G)
\rightarrow \prod_{v \in S} \textup{Br}_1(G_v) / (G(k_v)^{\bot})
\times \prod_{v \notin S \cup S'} \textup{Br}_1(G_v) /
(\mathcal{G}(\mathcal{O}_v)^{\bot})$$
Donc par dualit\'e locale, on en d\'eduit qu'il existe $(g'_v)_{v
  \notin S'} \in I := \prod_{v \in S} G(k_v)^{\wedge} \times \prod_{v
  \notin S \cup S'} \mathcal{G}(\mathcal{O}_v)$ tel que l'\'element $g
\in \prod_{v \in S} G(k_v)^{\wedge} \times \prod'_{v \notin S} G(k_v)$
obtenu en ``concat\'enant'' les $g_v$ et les $g'_v$ soit d'image nulle
dans $\textup{Br}_a(G)^D$. Or la suite 
$$1 \rightarrow \overline{G(k).\rho(G^{\textup{sc}}_{S_0})}^I
\rightarrow I \rightarrow \textup{Br}_a(G)^D$$
est exacte par le m\^eme argument que dans la preuve du corollaire
\ref{coro AF ouvert} (o\`u $\overline{(.)}^I$ d\'esigne l'adh\'erence
dans $I$). Par cons\'equent, $g'$ est dans
$\overline{G(k).\rho(G^{\textup{sc}}_{S_0})}^I$, et donc $(g_v)$ est
dans $\overline{H^0(U, \mathcal{G})}^{S'}$.
\end{dem}

\begin{rem}
{\rm
Dans ce lemme, on n'a en fait pas besoin que $\mathcal{G}$ soit r\'eductif sur $\textup{Spec}(\mathcal{O}_{k, S_0})$ :  ce r\'esultat est valable pour tout $\textup{Spec}(\mathcal{O}_{k, S_0})$-sch\'ema en groupes plat de type fini de fibre g\'en\'erique $G$.
}
\end{rem}

Poursuivons la preuve du th\'eor\`eme : soit $v$ une place de $S'$. On remarque d'abord le fait suivant :
\begin{lem}
Le groupe ab\'elien $H^0_{\textup{ab}}(\mathcal{O}_v, \mathcal{G})$ est $l$-divisible pour presque tout nombre premier $l$.
\end{lem}

\begin{dem}
On note $\mathcal{C}$ le complexe de sch\'emas en groupes commutatifs $[\mathcal{Z}^{\textup{sc}} \xrightarrow{\rho} \mathcal{Z}]$ associ\'e \`a $\mathcal{G}$, de sorte que $H^0_{\textup{ab}}(\mathcal{O}_v, \mathcal{G}) = \H^0(\mathcal{O}_v, \mathcal{C})$. On note aussi $\mathcal{K} := \textup{Ker}(\rho)$, qui est un sch\'ema en groupes commutatif fini, et $\mathcal{Z}' := \mathcal{Z} / \rho(\mathcal{Z}^{\textup{sc}})$.
\'Etant donn\'e un nombre premier $l$, on consid\`ere le diagramme commutatif de groupes ab\'eliens suivant :
\begin{displaymath}
\xymatrix{
H^1(\mathcal{O}_v, \mathcal{K}) \ar[r] \ar[d]^{[l]} & \H^0(\mathcal{O}_v, \mathcal{C}) \ar[r] \ar[d]^{[l]} & H^0(\mathcal{O}_v, \mathcal{Z}') \ar[r] \ar[d]^{[l]} & H^2(\mathcal{O}_v, \mathcal{K}) \ar[d]^{[l]} \\
H^1(\mathcal{O}_v, \mathcal{K}) \ar[r] & \H^0(\mathcal{O}_v, \mathcal{C}) \ar[r] & H^0(\mathcal{O}_v, \mathcal{Z}') \ar[r]  & H^2(\mathcal{O}_v, \mathcal{K})
}
\end{displaymath}
o\`u $[l]$ d\'esigne l'application de multiplication par $l$. La structure des groupes de Lie commutatifs $p$-adiques compacts assure que le morphisme $H^0(\mathcal{O}_v, \mathcal{Z}') \xrightarrow{[l]} H^0(\mathcal{O}_v, \mathcal{Z}')$ est surjectif pour presque tout $l$. Si $N$ d\'esigne le cardinal de $\mathcal{K}$, alors les groupes $\H^i(\mathcal{O}_v, \mathcal{K})$ sont de $N$-torsion, donc si $l$ est premier \`a $N$, le morphisme $H^i(\mathcal{O}_v, \mathcal{K}) \xrightarrow{[l]} H^i(\mathcal{O}_v, \mathcal{K})$ est un isomorphisme. Par cons\'equent, une chasse au diagramme assure imm\'ediatement que le morphisme $H^0(\mathcal{O}_v, \mathcal{C}) \xrightarrow{[l]} H^0(\mathcal{O}_v, \mathcal{C})$ est surjectif pour presque tout $l$, ce qui conclut la preuve.
\end{dem}

Sachant cela, on en d\'eduit qu'il suffit de trouver $S$ convenable tel que $\H^1_{S, S'}(k, \widehat{C})[l] = 0$ pour un ensemble fini donn\'e de nombres premiers $l$. Et pour cela, on utilise le th\'eor\`eme 7 de \cite{Ser}, section II.6.2, pour montrer par d\'evissage que $\H^1_{S_0, S'}(k, \widehat{C})[l]$ est fini pour chaque $l$, et on conclut comme dans la preuve du th\'eor\`eme 3 de \cite{HarAF} gr\^ace \`a l'hypoth\`ese sur la trivialit\'e de $\textup{\cyr{SH}}^1_{\omega}(k, \widehat{C}) / \textup{\cyr{SH}}^1(k, \widehat{C})$.

\end{dem}

\subsection{Une preuve ``g\'eom\'etrique'' \`a partir du cas des tores}
\label{section preuve 2}

Dans cette section, on va remontrer une partie du th\'eor\`eme \ref{theo fin} \`a partir du r\'esultat de Harari sur les 1-motifs (\cite{HarAF}, th\'eor\`eme 2). On va pour cela se ramener au cas connu des tores, en utilisant une z-extension.

Soit $k$ un corps de nombres, $G$ un $k$-groupe lin\'eaire r\'eductif tel que $(G^{\textup{sc}})^i_{S_0}$ est non compact pour tout facteur presque $k$-simple $(G^{\textup{sc}})^i$ de $G^{\textup{sc}}$. On sait que $G$ admet une z-extension (voir \cite{MS}, p 297, proposition 3.1), c'est-\`a-dire qu'il existe une suite exacte
$$1 \rightarrow P \rightarrow H \rightarrow G \rightarrow 1$$
de groupes r\'eductifs sur $k$, de sorte que $P$ est un $k$-tore quasi-trivial et $H$ un $k$-groupe r\'eductif tel que son sous-groupe d\'eriv\'e $H^{\textup{ss}}$ soit (semi-simple) simplement connexe.
On consid\`ere alors le diagramme commutatif suivant, dont les lignes sont exactes :
\begin{displaymath}
\xymatrix{
H^0(k, P) \ar[r] \ar[d] & H^0(k, H) \ar[r] \ar[d] & H^0(k, G) \ar[r] \ar[d] & H^1(k, P) = 1 \\
P^0(k, P) \ar[r] \ar[d] & P^0(k, H) \ar[r] \ar[d] & P^0(k, G) \ar[r] \ar[d] & P^1(k, P) = 1 \\
\textup{Br}_a(P)^D = H^2(k, \widehat{P})^D \ar[r] & \textup{Br}_a(H)^D \ar[r] & \textup{Br}_a(G)^D & 
}
\end{displaymath}
Le point crucial est le suivant : le morphisme $P^0(k, P) \rightarrow
H^2(k, \widehat{P})^D$ est surjectif. En effet, $P$ \'etant
quasi-trivial, il suffit de montrer la surjectivit\'e de $P^0(L, \G_m)
\rightarrow H^2(L, \Z)^D$ pour $L$ un corps de nombres, et ceci est
une cons\'equence de la th\'eorie du corps de classes : en effet, cette
fl\`eche s'identifie, gr\^ace \`a la proposition 3 du paragraphe 1,
chapitre XIV de \cite{Ser2}, au morphisme de r\'eciprocit\'e de la
th\'eorie du corps de classes $\A_L^* \rightarrow
\Gamma_L^{\textup{ab}}$ entre le groupe des classes d'id\`eles et
l'ab\'elianis\'e du groupe de Galois absolu de $L$; et il est bien
connu que cette fl\`eche est surjective (voir par exemple \cite{NSW}, proposition 8.1.25).

Supposons d\'esormais le r\'esultat connu pour le groupe $H$, i.e. $\overline{H(k).H^{\textup{sc}}_{S_0}} = P^0(k, H)^{\textup{Br}_a}$ avec les notations pr\'ec\'edentes.
D\'eduisons-en le r\'esultat analogue pour $G$ : prenons $(g_v) \in P^0(k, G)^{\textup{Br}_a}$. Par exactitude du diagramme pr\'ec\'edent, on rel\`eve $(g_v)$ en $(h_v) \in P^0(k, H)$. A priori, $(h_v)$ n'est pas dans $P^0(k, H)^{\textup{Br}_a}$, mais on va le modifier par un \'el\'ement de $P^0(k, P)$ pour que cela soit le cas. En effet, l'image de $(h_v)$ dans $\textup{Br}_a(H)^D$ s'envoie sur $0$ dans $\textup{Br}_a(G)^D$, donc elle se rel\`eve en un \'el\'ement $\delta$ dans $\textup{Br}_a(P)^D$. Par surjectivit\'e de la fl\`eche $P^0(k, P) \rightarrow \textup{Br}_a(P)^D$, $\delta$ se rel\`eve en $(p_v) \in P^0(k, P)$. Notons alors $(h'_v) := (p_v)^{-1}.(h_v) \in P^0(k, H)$. Par construction, $(h'_v)$ s'envoie sur $0$ dans $\textup{Br}_a(H)^D$, donc il est dans l'adh\'erence de $H(k).H^{\textup{sc}}_{S_0}$. Et $(h'_v)$ rel\`eve $(g_v) \in P^0(k, G)$, donc $(g_v) \in \overline{G(k).G^{\textup{sc}}_{S_0}}$.

Il suffit donc d\'esormais de montrer le r\'esultat pour un groupe r\'eductif $H$ tel que $H^{\textup{ss}}$ est simplement connexe \footnote{Je remercie Jean-Louis Colliot-Th\'el\`ene pour ses suggestions \`a propos de cette partie de la preuve.}. Pour cela, on consid\`ere la suite exacte
$$1 \rightarrow H^{\textup{ss}} \rightarrow H \rightarrow T \rightarrow 1$$
o\`u $T$ est un $k$-tore. 0n voit donc $T$ comme un espace homog\`ene de $H$ \`a stabilisateur $H^{\textup{ss}}$.

On se donne un point $(h_v) \in P^0(k, H)$ Brauer-Manin orthogonal au groupe $\textup{Br}_a(H)$, que l'on pousse dans $P^0(k, T)$ pour obtenir un point $(t_v)$. Gr\^ace au r\'esultat de Harari sur les tores (th\'eor\`eme 2 de \cite{HarAF}), on sait que $(t_v)$ est dans l'adh\'erence forte de $T(k)$. Soit alors un point $t_0$ suffisamment proche de $(t_v)$ pour la topologie forte sur $P^0(k, T)$. La fibre $H_{t_0}$ est alors un espace principal homog\`ene de $H^{\textup{ss}}$ (qui est simplement connexe), et cette fibre a des points locaux en toutes les places r\'eelles par approximation, donc elle a un point rationnel puisque $H^1(k, G^{\textup{ss}}) \rightarrow \prod_{v \textup{ r\'eelle}} H^1(k_v, G^{\textup{ss}})$ est injective (principe de Hasse pour les groupes semi-simples simplement connexes, voir \cite{PR}, th\'eor\`eme 6.6 par exemple), donc $H_{t_0}$ est $k$-isomorphe \`a un groupe semi-simple simplement connexe. On conclut alors par le th\'eor\`eme d'approximation forte sur $H_{t_0}$ (et le th\'eor\`eme des fonctions implicites pour obtenir des points locaux dans $H_{t_0}$ proches des points $h_v$ initiaux).

Finalement, on a retrouv\'e le fait que $\overline{G(k).G^{\textup{sc}}_{S_0}} = P^0(k, G)^{\textup{Br}_a}$ pour tout groupe r\'eductif $G$, \`a partir du r\'esultat analogue pour les tores.

\section{Groupe des classes d'un groupe r\'eductif}
\label{section nbre classes}
Dans cette section, le probl\`eme est le suivant :  \'etant donn\'e un
$k$-groupe r\'eductif, et $S_0$ un ensemble fini de places contenant
les places archim\'ediennes tel que $G^{\textup{sc}}$ v\'erifie
l'approximation forte hors de $S_0$, et si $\mathcal{G}$ est un
sch\'ema en groupes plat de type fini sur $\mathcal{O}_k$,
de fibre g\'en\'erique $G$, on d\'efinit l'ensemble des classes de
$\mathcal{G}$ comme l'ensemble de doubles classes suivant :
$$\textup{Cl}_{S_0}(\mathcal{G}) := G(k) \backslash G(\A_k) / \mathcal{G}(\A(S_0))$$
o\`u $ \mathcal{G}(\A(S_0)) := \prod_{v \in S_0} G(k_v) \times
\prod_{v \notin S_0} \mathcal{G}(\mathcal{O}_v)$, et on note $\textup{cl}_{S_0}(\mathcal{G})$ son cardinal.

\begin{rem}
{\rm Le mod\`ele $\mathcal{G}$ n'est pas n\'ecessairement s\'epar\'e : pour toute place
  $v$ de $k$, on dispose d'un morphisme naturel
  $\mathcal{G}(\mathcal{O}_v) \rightarrow G(k_v)$ qui n'est pas
  n\'ecessairement injectif. Aussi le quotient $G(\A_k) / \mathcal{G}(\A(S_0))$ d\'esigne-t-il le quotient de
  $G(\A_k)$ par l'image de $\mathcal{G}(\A(S_0))$ dans $G(\A_k)$.}
\end{rem}

On sait que $\textup{Cl}_{S_0}(\mathcal{G})$ est naturellement un groupe ab\'elien fini (voir \cite{PR}, proposition 8.8 dans le cas o\`u $G$ est semi-simple, ou \cite{Kne} pour le cas g\'en\'eral), de cardinal $\textup{cl}_{S_0}(\mathcal{G})$, et qu'il co\"incide avec le groupe $G(\A_k) / G(k).\mathcal{G}(\A(S_0))$.
L'objectif ici est de d\'ecrire ce groupe "explicitement" en fonction du groupe de Brauer de $G$.

Pour cela, on d\'efinit le groupe $\textup{Br}^{S_0}_{a,
  \textup{nr}}(\mathcal{G})$ comme \'etant le sous-groupe des
\'el\'ements de $\textup{Br}_a(G)$ (identifi\'e \`a $\textup{Br}_{e}(G)
:= \textup{Ker}( \textup{Br}_1(G) \xrightarrow{e^*} \textup{Br}(k))$,
$e \in G(g)$ \'etant l'\'el\'ement neutre de $G$) qui sont orthogonaux \`a $G(k_v)$ dans $\textup{Br}_1(G_v)$ pour $v \in S_0$, et qui sont orthogonaux \`a $\mathcal{G}(\mathcal{O}_v)$ pour $v \notin S_0$, pour l'accouplement local $G(k_v) \times \textup{Br}_1(G_v) \rightarrow \Q / \Z$.

\begin{theo}
\label{theo classes}
L'accouplement de Brauer-Manin induit une dualit\'e parfaite de groupes finis, fonctorielle en $\mathcal{G}$ :
$$\textup{Cl}_{S_0}(\mathcal{G}) \times \textup{Br}^{S_0}_{a, \textup{nr}}(\mathcal{G}) / \textup{\cyr{B}}(G) \rightarrow \Q / \Z$$
\end{theo}

\begin{dem} 
Si $U$ est un ouvert suffisamment petit, la suite suivante :
$$\overline{H^0(U, \mathcal{G}).G^{\textup{sc}}_{S_0}} .  \mathcal{P}^0_{S_0}(\mathcal{G}) \rightarrow \mathcal{P}^0_S(\mathcal{G}) \rightarrow \textup{Br}^{S_0}_{a, \textup{nr}}(\mathcal{G})^D$$
est exacte : en effet, cela se d\'eduit facilement de la suite exacte du corollaire \ref{coro AF ouvert} et de la dualit\'e locale (voir la preuve du lemme \ref{lem approx faible}).
On passe ensuite \`a la limite sur $U$. On obtient la suite exacte suivante :
$$\overline{G(k).G^{\textup{sc}}_{S_0}} .\mathcal{P}^0_{S_0}(\mathcal{G}) \rightarrow P^0(k, G) \rightarrow \textup{Br}^{S_0}_{a, \textup{nr}}(\mathcal{G})^D$$
Or $\mathcal{P}^0_{S_0}(\mathcal{G})$ est un sous-groupe ouvert pour la topologie ad\'elique, on en d\'eduit donc facilement que $\overline{G(k).G^{\textup{sc}}_{S_0}} .\mathcal{P}^0_{S_0}(\mathcal{G})$ co\"incide avec l'image de $G(k) . \mathcal{P}^0_{S_0}(\mathcal{G})$ dans $P^0(k, G)$, d'o\`u une suite exacte
$$G(k) . \mathcal{P}^0_{S_0}(\mathcal{G}) \rightarrow P^0(k, G) \rightarrow \textup{Br}^{S_0}_{a, \textup{nr}}(\mathcal{G})^D$$
Pour finir, il est clair que l'on peut remplacer simultan\'ement $ \mathcal{P}^0_{S_0}(\mathcal{G})$ par $\mathcal{G}(\A(S_0))$ et $P^0(k, G)$ par $G(\A_k)$, d'o\`u la suite exacte
$$G(k) . \mathcal{G}(\A(S_0)) \rightarrow G(\A_k) \rightarrow \textup{Br}^{S_0}_{a, \textup{nr}}(\mathcal{G})^D$$
et enfin on v\'erifie imm\'ediatement (gr\^ace au th\'eor\`eme
\ref{theo fin}) que le conoyau de la derni\`ere fl\`eche s'identifie \`a $\textup{\cyr{B}}(G)^D$, ce qui conclut la preuve du th\'eor\`eme \ref{theo classes}.
\end{dem}

Une cons\'equence de ce r\'esultat est une g\'en\'eralisation et un raffinement du th\'eor\`eme 8.12 de \cite{PR}, qui concerne la comparaison du nombre de classes d'un groupe avec celui de ses tores maximaux :

\begin{cor}
Si $T$ est un $k$-tore maximal de $G$, et si $\mathcal{T}$ est l'adh\'erence sch\'ematique de $T$ dans $\mathcal{G}$, on a l'in\'egalit\'e suivante :
$$\textup{cl}_{S_0}(\mathcal{T}) \geq \frac{| \textup{\cyr{SH}}^1(G)|}{| \textup{\cyr{SH}}^1(T)| . | H^1(k, \widehat{T^{\textup{sc}}})|}  \textup{cl}_{S_0}(\mathcal{G})$$
o\`u $H^1(k, \widehat{T^{\textup{sc}}})$ est isomorphe \`a $\textup{Pic}(T^{\textup{sc}})$. On a \'egalement l'in\'egalit\'e :
$$\textup{cl}_{S_0}(\mathcal{T}) \geq \frac{| \textup{Pic}(T) / \textup{Pic}(G) |}{| \textup{Pic}(T^{\textup{sc}}) | }  \frac{| \textup{\cyr{SH}}^1(G)|}{| \textup{\cyr{SH}}^1(T)|}  \textup{cl}_{S_0}(\mathcal{G})$$
De fa\c con plus pr\'ecise, il existe un entier $N(\mathcal{G}) \geq 1$ tel que 
$$\textup{cl}_{S_0}(\mathcal{T}) = N(\mathcal{G}) .  \frac{| \textup{\cyr{SH}}^1(G)|}{| \textup{\cyr{SH}}^1(T)| . | H^1(k, \widehat{T^{\textup{sc}}})|}  \textup{cl}_{S_0}(\mathcal{G})$$
\end{cor}

\begin{dem}
Il s'agit de comparer les cardinaux des ensembles $H^2_{\textup{nr}, S_0}(k, \widehat{\mathcal{T}})$ et $\H^1_{\textup{nr}, S_0}(k, \widehat{\mathcal{C}})$, o\`u $H^2_{\textup{nr}, S_0}(k, \widehat{\mathcal{T}})$ est le sous-groupe de $H^2(k, \widehat{T})$ form\'e des \'el\'ements $\alpha$ tels que $\alpha_v = 0$ si $v \in S_0$ et $\alpha_v$ est orthogonal \`a $\mathcal{T}(\mathcal{O}_v)$ pour $v \notin S_0$, et $\H^1_{\textup{nr}, S_0}(k, \widehat{\mathcal{C}})$ est le sous-groupe de $\H^1(k, \widehat{C})$ form\'e des \'el\'ements $\beta$ tels que $\beta_v$ est orthogonal \`a $G(k_v)$ si $v \in S_0$ et \`a $\mathcal{G}(\mathcal{O}_v)$ pour $v \notin S_0$.
Pour cela, on remarque que la suite exacte
$$H^1(k, \widehat{T^{\textup{sc}}}) \rightarrow \H^1(k, \widehat{C}) \rightarrow H^2(k, \widehat{T})$$
envoie le sous-groupe $\H^1_{\textup{nr}, S_0}(k, \widehat{\mathcal{C}})$ dans $H^2_{\textup{nr}, S_0}(k, \widehat{\mathcal{T}})$ en raison de la commutativit\'e des diagrammes suivants issus de la dualit\'e locale :
\begin{displaymath}
\xymatrix{
G(k_v) \times \H^1(k_v, \widehat{C}) \ar@<5ex>[d] \ar[r] & \Q / \Z \ar[d]^{=} \\
T(k_v) \times H^2(k_v, \widehat{T}) \ar[r] \ar@<7ex>[u] & \Q / \Z 
}
\end{displaymath}
On en d\'eduit donc l'in\'egalit\'e suivante sur les cardinaux de ces ensembles :
$$| \H^1_{\textup{nr}, S_0}(k, \widehat{\mathcal{C}}) | \leq | H^2_{\textup{nr}, S_0}(k, \widehat{\mathcal{T}}) | . | H^1(k, \widehat{T^{\textup{sc}}}) |$$
et on conclut gr\^ace au th\'eor\`eme \ref{theo classes}, apr\`es avoir identifi\'e $\H^1_{\textup{nr}, S_0}(k, \widehat{\mathcal{C}})  \cong \textup{Br}^{S_0}_{a, \textup{nr}}(\mathcal{G}) $ et $ H^2_{\textup{nr}, S_0}(k, \widehat{\mathcal{T}}) \cong  \textup{Br}^{S_0}_{a, \textup{nr}}(\mathcal{T})$, ce qui donne la premi\`ere in\'egalit\'e. Pour la seconde, on consid\`ere la suite exacte suivante :
$$\H^0(k, \widehat{C}) \rightarrow H^1(k, \widehat{T}) \rightarrow H^1(k, \widehat{T^{\textup{sc}}}) \rightarrow \H^1_{S_0}(k, \widehat{\mathcal{C}})' \rightarrow H^2_{\textup{nr},S_0}(k, \widehat{\mathcal{T}})$$
o\`u $\H^1_{S_0}(k, \widehat{\mathcal{C}})'$ d\'esigne le sous-groupe de $\H^1(k, \widehat{C})$, contenant $\H^1_{S_0}(k, \widehat{\mathcal{C}})$, form\'e des \'el\'ements dans $\H^1(k, \widehat{C})$ localement orthogonaux \`a l'image de $T(k_v)$ (ou $\mathcal{T}(\mathcal{O}_v)$ selon que $v \in S_0$ ou non). Cette suite exacte fournit l'in\'egalit\'e suivante sur les cardinaux :
$$| H^1(k, \widehat{T^{\textup{sc}}}) | | H^2_{\textup{nr}, S_0}(k, \widehat{\mathcal{T}}) | \geq | H^1(k, \widehat{T}) / \H^0(k, \widehat{C}) | |\H^1_{S_0}(k, \widehat{\mathcal{C}})' |$$
et on conclut en identifiant $H^1(k, \widehat{T^{\textup{sc}}})$,
$H^1(k, \widehat{T})$ et $\H^0(k, \widehat{C})$ \`a
$\textup{Pic}(T^{\textup{sc}})$, $\textup{Pic}(T)$ et $\textup{Pic}(G)$
respectivement, en minorant $|\H^1_{S_0}(k, \widehat{\mathcal{C}})' |$
par $|\H^1_{\textup{nr},S_0}(k, \widehat{\mathcal{C}}) |$, et en
appliquant le th\'eor\`eme \ref{theo classes} pour exprimer
$|\H^1_{\textup{nr}, S_0}(k, \widehat{\mathcal{C}}) |$ et $ |
H^2_{\textup{nr}, S_0}(k, \widehat{\mathcal{T}}) | $ en fonction des nombres de classes respectifs de $T$ et $G$.
\end{dem}

Une autre application concerne la comparaison du nombre de classes de $G$ avec celui d'un quotient de $G$ :

\begin{cor} Soit $1 \rightarrow G_1 \rightarrow G_2 \rightarrow G_3
  \rightarrow 1$ une suite exacte de $k$-groupes r\'eductifs. Soit
  $\mathcal{G}_2$ un sch\'ema en groupes plat de type fini
  \'etendant $G_2$ sur $\textup{Spec}(\mathcal{O}_{k, S_0})$. On note
  $\mathcal{G}_1$ l'adh\'erence sch\'ematique de $G_1$ dans
  $\mathcal{G}_2$ et $\mathcal{G}_3$ le quotient de $\mathcal{G}_2$
  par $\mathcal{G}_1$. Alors $\mathcal{G}_3$ est un sch\'ema en
  groupes plat de type fini, et on a les in\'egalit\'es suivantes :
$$\textup{cl}_{S_0}(\mathcal{G}_3) \leq \frac{| \textup{\cyr{SH}}^1(G_2)|}{| \textup{\cyr{SH}}^1(G_3)|} | \textup{Pic}(G_1) |  \textup{cl}_{S_0}(\mathcal{G}_2)$$
et
$$\textup{cl}_{S_0}(\mathcal{G}_3) \leq \frac{| \textup{\cyr{SH}}^1(G_2)|}{| \textup{\cyr{SH}}^1(G_3)|} \frac{| \textup{Pic}(G_1) |}{| \textup{Pic}(G_2) / \textup{Pic}(G_3) |}  \textup{cl}_{S_0}(\mathcal{G}_2)$$
Plus pr\'ecis\'ement, il existe un entier $M(\mathcal{G}_2) \geq 1$ tel que 
$$ \textup{cl}_{S_0}(\mathcal{G}_2) = M(\mathcal{G}_2) . \frac{| \textup{\cyr{SH}}^1(G_3)|}{| \textup{\cyr{SH}}^1(G_2)|} \frac{1}{ | \textup{Pic}(G_1) | }  \textup{cl}_{S_0}(\mathcal{G}_3)$$
\end{cor}

\begin{dem}
On sait que $\mathcal{G}_1$ est un sous-sch\'ema en groupes plat
de type fini distingu\'e dans $\mathcal{G}_2$. Alors par
\cite{SGA3}, expos\'e $\textup{VI}_{\textup{B}}$, proposition 9.2, et
gr\^ace au th\'eor\`eme 4.C. de \cite{Ana}, le quotient $\mathcal{G}_2
/ \mathcal{G}_1$ est repr\'esentable par un
$\textup{Spec}(\mathcal{O}_{k, S_0})$-sch\'ema en groupes plat de type fini $\mathcal{G}_3$.

On consid\`ere alors la suite exacte suivante (voir par exemple \cite{San}, corollaire 6.11) :
$$\textup{Pic}(G_3) \rightarrow \textup{Pic}(G_2) \rightarrow \textup{Pic}(G_1) \rightarrow \textup{Br}_a(G_3) \rightarrow \textup{Br}_a(G_2) \rightarrow \textup{Br}_a(G_1)$$
On en d\'eduit alors la suite exacte suivante :
$$\textup{Pic}(G_3) \rightarrow \textup{Pic}(G_2) \rightarrow \textup{Pic}(G_1)' \rightarrow \textup{Br}^{S_0}_{a, \textup{nr}}(\mathcal{G}_3) \rightarrow \textup{Br}^{S_0}_{a, \textup{nr}}(\mathcal{G}_2)$$
o\`u $\textup{Pic}(G_1)'$ est un sous-groupe de $\textup{Pic}(G_1)$.
Le corollaire est alors une cons\'equence directe du th\'eor\`eme \ref{theo classes}, en exprimant les cardinaux des groupes apparaissant dans cette suite exacte en fonction des nombres de classes respectifs.
\end{dem}

\begin{exs}
Le corollaire pr\'ec\'edent fournit entre autres les cas particuliers suivants :
\begin{itemize}
	\item Si $G$ est un $k$-groupe r\'eductif, et si $\mathcal{G}$
          est un sch\'ema en groupes plat de type fini
          \'etendant $G$ sur $\textup{Spec}(\mathcal{O}_{k, S_0})$, on
          note $\mathcal{G}^{\textup{tor}}$ le quotient de
          $\mathcal{G}$ par l'adh\'erence de $G^{\textup{ss}}$. C'est
          un sch\'ema en groupes plat de type fini sur $\textup{Spec}(\mathcal{O}_{k, S_0})$, dont la fibre g\'en\'erique est le tore $G^{\textup{tor}} := G / G^{\textup{ss}}$. Alors on a l'in\'egalit\'e 
$$\textup{cl}_{S_0}(\mathcal{G}^{\textup{tor}}) \leq \frac{| \textup{\cyr{SH}}^1(G)|}{| \textup{\cyr{SH}}^1(G^{\textup{tor}})|} | \textup{Pic}(G^{\textup{ss}}) |  \textup{cl}_{S_0}(\mathcal{G})$$
et m\^eme
$$\textup{cl}_{S_0}(\mathcal{G}) = M(\mathcal{G}) . \frac{| \textup{\cyr{SH}}^1(G^{\textup{tor}})|}{| \textup{\cyr{SH}}^1(G)|} \frac{1}{| \textup{Pic}(G^{\textup{ss}}) |} \textup{cl}_{S_0}(\mathcal{G}^{\textup{tor}})$$ pour un entier $M(\mathcal{G}) \geq 1$.
	\item Si $G$ est un groupe r\'eductif, tel que $G^{\textup{ss}}$ est simplement connexe, alors avec les notations du point pr\'ec\'edent, $\textup{cl}_{S_0}(\mathcal{G}^{\textup{tor}})$ divise $\textup{cl}_{S_0}(\mathcal{G})$, et on dispose m\^eme d'un morphisme canonique surjectif
$$\textup{Cl}_{S_0}(\mathcal{G}) \rightarrow \textup{Cl}_{S_0}(\mathcal{G}^{\textup{tor}})$$
	\item Si $1 \rightarrow Z \rightarrow H \rightarrow G \rightarrow 1$ est une z-extension de $G$, et si $\mathcal{H}$ est un sch\'ema en groupes plat de type fini \'etendant $H$ sur $\textup{Spec}(\mathcal{O}_{k, S_0})$, alors on a un morphisme canonique surjectif
$$\textup{Cl}_{S_0}(\mathcal{H}) \rightarrow \textup{Cl}_{S_0}(\mathcal{G})$$
\end{itemize}
\end{exs}

Un cas particulier du th\'eor\`eme \ref{theo classes} est celui des tores :

\begin{cor}
Soit $T/k$ un tore alg\'ebrique. Soit $S_0$ un ensemble fini de places contentant les places archim\'ediennes. Alors on a une dualit\'e parfaite de groupes finis
$$\textup{Cl}_{S_0}(\mathcal{T}) \times H^2_{\textup{nr}, S_0}(k, \widehat{\mathcal{T}}) / \textup{\cyr{SH}}^2(\widehat{T}) \rightarrow \Q / \Z$$
o\`u $H^2_{\textup{nr}, S_0}(k, \widehat{\mathcal{T}})$ est le sous-groupe de $H^2(k, \widehat{T})$ form\'e des \'el\'ements $\alpha$ tels que $\alpha_v = 0$ si $v \in S_0$ et $\alpha_v$ est orthogonal \`a $\mathcal{T}(\mathcal{O}_v)$ pour $v \notin S_0$.
En particulier, si $\mathcal{T}$ a bonne r\'eduction hors de $S_0$, on obtient
$$H^2_{\textup{nr}, S_0}(k, \widehat{\mathcal{T}}) = \left\{ \alpha \in H^2(k, \widehat{T}) : \alpha_v = 0 \textup{ si }v \in S_0 \textup{ et }  \alpha_v \textup{ non ramifi\'e si }v \notin S_0 \right\}$$
\end{cor}

\begin{ex}
Un calcul imm\'ediat \`a l'aide du corps de classes de Hilbert montre
que ce corollaire est coh\'erent avec le fait bien connu que $\textup{Cl}_{S_{\infty}}(\G_m)$ est
isomorphe au groupe des classes d'id\'eaux de $k$.
\end{ex}

Un autre cas particulier du th\'eor\`eme \ref{theo classes} est le suivant :

\begin{cor}
Soit $G$ un $k$-groupe semi-simple de groupe fondamental $B$. Soit
$\mathcal{G}$ un $\textup{Spec}(\mathcal{O}_k)$-sch\'ema en groupes
plat de type fini \'etendant $G$. Alors on a une dualit\'e parfaite
$$\textup{Cl}_{S_0}(\mathcal{G}) \times H^1_{\textup{nr}, S_0}(k, \widehat{\mathcal{B}}) / \textup{\cyr{SH}}^1(\widehat{B}) \rightarrow \Q / \Z$$
o\`u $H^1_{\textup{nr}, S_0}(k, \widehat{\mathcal{B}})$ est l'ensemble des $\alpha \in H^1(k, \widehat{B})$ tels que $\alpha_v$ est orthogonal \`a l'image de $G(k_v)$ dans $H^1(k_v, B)$ si $v \in S_0$, et $\alpha_v$ est orthogonal \`a l'image de $\mathcal{G}(\mathcal{O}_v)$ si $v \notin S_0$.
\end{cor}

Terminons cette partie par une remarque sur le morphisme de normes d\'efini dans \cite{Tharec} (voir th\'eor\`eme 14 de \cite{Tharec}). \`A l'aide des r\'esultats pr\'ec\'edents, on peut montrer le th\'eor\`eme suivant :

\begin{theo}[Thang]
Supposons que $\mathcal{G}$ soit un sch\'ema en groupes r\'eductifs sur $\mathcal{O}_{k, S_0}$. Alors pour toute extension finie $L / k$, il existe un homomorphisme canonique et fonctoriel en $\mathcal{G}$ :
$$N_{L / k} : \textup{Cl}_{S_{0,L}}(\mathcal{G}_{\mathcal{O}_L}) \rightarrow \textup{Cl}_{S_0}(\mathcal{G})$$
de sorte que pour une tour d'extensions $k \subset L \subset M$, on ait 
$$N_{M / k} = N_{L / k} \circ N_{M / L}$$
Ce morphisme $N_{L / k}$ est induit par le dual du morphisme de restriction usuel sur le groupe de Brauer $\textup{Res}_{L / k} : \textup{Br}(G) \rightarrow \textup{Br}(G_L)$, via le th\'eor\`eme \ref{theo classes}.
\end{theo}

\begin{dem}
En regard du th\'eor\`eme \ref{theo classes}, il suffit de montrer que le morphisme de restriction $\textup{Res}_{L / k} : \textup{Br}(G) \rightarrow \textup{Br}(G_L)$ envoie bien le sous-groupe $\textup{Br}_{a, \textup{nr}}^{S_0}(\mathcal{G})$ de $\textup{Br}(G)$ dans le sous-groupe $\textup{Br}_{a, \textup{nr}}^{S_{0,L}}(\mathcal{G}_{\mathcal{O}_L})$ de $\textup{Br}(G_L)$. Et ce fait est une cons\'equence imm\'ediate des deux r\'esultats suivants :
\begin{itemize}
	\item l'existence d'un morphisme de corestriction $H^1_{\textup{ab}}(\mathcal{O}_w, \mathcal{G}) \rightarrow H^1_{\textup{ab}}(\mathcal{O}_v, \mathcal{G})$ pour toute place $w$ de $L$ divisant une place $v$ de $k$ hors de $S_0$ (voir par exemple \cite{Del} et \cite{SGA4}, expos\'e XVII, section 6.3).
	\item la validit\'e du principe de corestriction pour le morphisme d'ab\'elianisation $\textup{ab}^0 : G(k_v) \rightarrow H^0_{\textup{ab}}(k_v, G)$ pour $v \in S_0$ (voir par exemple \cite{Del} ou \cite{Tha}).
\end{itemize}
\end{dem}

\section{Une suite de Poitou-Tate non ab\'elienne}
\label{section PTNA}
On se propose dans cette section de d\'emontrer l'existence d'une
suite de Poitou-Tate pour un groupe lin\'eaire connexe $G / k$,
fonctorielle en $G$, et qui est compatible, via les fl\`eches
d'ab\'elianisation, \`a la suite exacte de Poitou-Tate pour
le complexe de $k$-tores $C_G = [T_G^{\textup{sc}} \rightarrow T_G]$
(voir th\'eor\`eme 6.1 de \cite{Dem1}).

Concernant la d\'efinition et les
propri\'et\'es de base de l'ensemble de cohomologie non ab\'elienne
$H^2(k, G)$, on renvoie \`a \cite{FSS}, section 1, ou alors \`a
\cite{Bor2}. On rappelle que l'ensemble $H^2(k, G)$, lorsqu'il est
non vide (ce qui est le cas ici), est un espace
principal homog\`ene sous le groupe $H^2(k, Z(G))$ (pour une action
naturelle, not\'ee $+$), $Z(G)$ \'etant le
centre de $G$. Ainsi, si $\eta, \eta' \in H^2(k, G)$, il existe un
unique $z \in H^2(k, Z)$ tel que $z + \eta' = \eta$; on note alors
$z := \eta - \eta'$. On consid\`ere alors le complexe de $k$-groupes
commutatifs $C := [Z^{\textup{sc}} \rightarrow Z]$, quasi-isomorphe \`a
$[T^{\textup{sc}} \rightarrow T]$. On a un morphisme naturel $j_* :
H^2(k, Z) \rightarrow \H^2(k, C)$, et si $n(G) \in H^2(k, G)$ est la
classe neutre associ\'ee \`a la $k$-forme $G$, on a par d\'efinition
$\textup{ab}^2(\eta) := j_*(\eta - n(G)) \in H^2_{\textup{ab}}(k, G)
:= \H^2(k, C)$ (voir \cite{Bor2}, section 5.3).

D\'efinissons une relation d'\'equivalence $\sim$ sur le $H^2$
non ab\'elien d'un $k$-groupe connexe $G$ de centre $Z$ : soient
$\eta, \eta' \in H^2(k, G)$. On dira que $\eta$ et $\eta'$ sont
\'equivalents lorsque $\eta - \eta' \in H^2(k, Z)$ est dans l'image du
cobord $\delta : H^1(k, G / Z) \rightarrow H^2(k, Z)$.
La classe d\'equivalence de la classe neutre $n(G)$ associ\'ee \`a la
$k$-forme $G$ est donc exactement l'ensemble des classes neutres de $H^2(k,
G)$ (voir \cite{Bor2}, th\'eor\`eme 5.5). On v\'erifie alors (voir \cite{Bor2}, section 5) que l'on a la
relation suivante : $\eta \sim \eta'$ \SSI $\textup{ab}^2(\eta) =
\textup{ab}^2(\eta') \in H^2_{\textup{ab}}(k, G)$, donc il est clair que
$\sim$ est une relation d'\'equivalence. 

Enfin, on rappelle que $\textup{UPic}(\overline{G})$ est d\'efini comme le complexe de modules galoisiens (en degr\'es $-1$ et $0$)
$$\textup{UPic}(\overline{G}) := [\overline{k}(G)^* / \overline{k}^* \rightarrow \textup{Div}(\overline{G})]$$
qui est quasi-isomorphe au dual du complexe de $k$-tores $[T_G^{\textup{sc}} \rightarrow T_G]$ par le corollaire 2.20 de \cite{BVH} ($T_G$ et $T_G^{\textup{sc}}$ sont des tores maximaux de $G^{\textup{red}}$ et $G^{\textup{sc}}$). On obtient alors le r\'esultat suivant :
\begin{theo}
\label{PTNA}
Soit $G / k$ un groupe lin\'eaire connexe.
\begin{itemize}
	\item Alors on a une suite exacte naturelle d'ensembles
          point\'es (la premi\`ere ligne est une suite exacte de groupes), fonctorielle en $G$ :
\begin{displaymath}
\xymatrix{
& (\textup{Br}_a G )^D \ar[d] & P^0(k, G) \ar[l] & P^0(k, G)^{\textup{Br}_a} \ar[l] & 0 \ar[l] \\
& H^1(k, G) \ar[r] & P^1(k, G) \ar[r] & (\textup{Pic } G )^D \ar[d] & \\
0 & (k[G]^* / k^*)^D \ar[l] & \bigoplus_{v \in \Omega_k} H^2(k_v, G) / \sim \ar[l] & H^2(k, G) / \sim \ar[l] &
}
\end{displaymath}
et cette suite est compatible avec la suite exacte de Poitou-Tate pour
le complexe de $k$-tores $C_G = [T_G^{\textup{sc}} \rightarrow T_{G}]$ du
th\'eor\`eme 6.1 de \cite{Dem1}, via les applications
d'ab\'elianisation $\textup{ab}^i_G : H^i(.,G) \rightarrow \H^i(.,C_G)$.
	\item Si de plus $G$ v\'erifie que $G^i_{S_0} := \prod_{v \in S_0} G^i(k_v)$  est non compact pour tout $k$-facteur presque $k$-simple $G^i$ de $G^{\textup{sc}}$, alors on peut identifier $ P^0(k, G)^{\textup{Br}_a}$ \`a $\overline{G^{\textup{scu}}_{S_0} G(k)}$.
	\item Dans tous les cas, on dispose de la suite exacte duale de groupes ab\'eliens :
\begin{displaymath}
\xymatrix{
0 \ar[r] & \widehat{G}(k)^{\wedge} \ar[r] & \left( \prod' \widehat{G}(k_v) \right)^{\wedge} \ar[r] & H^2_{\textup{ab}}(k, G)^D \ar[d] \\
& H^1_{\textup{ab}}(k, G)^D \ar[d] & \prod_{v \in \Omega_k} \textup{Pic}(G_v) \ar[l] & \textup{Pic}(G) \ar[l] \\
& \textup{Br}_a(G) \ar[r] & \left( \prod'_v \textup{Br}_a(G_v)
\right)_{\textup{tors}} \ar[r] & \left( H^0_{\textup{ab}}(k, G)^D \right)_{\textup{tors}} \ar[d] \\
0 & \textup{Ker}(\rho)(k)^D \ar[l] & \P^2(k,\textup{UPic}(\overline{G})) \ar[l] & \H^2(k, \textup{UPic}(\overline{G})) \ar[l]
}
\end{displaymath}
o\`u le produit restreint $\prod'_v \textup{Br}_a(G_v)$ est
consid\'er\'e par rapport aux sous-groupes $\H^1(\mathcal{O}_v,
[\widehat{\mathcal{T}_G} \rightarrow
\widehat{\mathcal{T}_G^{\textup{sc}}}])$.
\end{itemize}
\end{theo}

\begin{exs}
\begin{itemize}
	\item Si $G = T$ est un $k$-tore, on obtient la suite de
          Poitou-Tate "usuelle" pour les tores (voir par exemple
          \cite{HSz}, th\'eor\`eme 5.6 ou \cite{Dem1}, th\'eor\`eme 6.1) :
\begin{displaymath}
\xymatrix{
& H^2(k, \widehat{T})^D \ar[d] & P^0(k, T) \ar[l] & \overline{T(k)} \ar[l] & 0 \ar[l] \\
& H^1(k, T) \ar[r] & P^1(k, T) \ar[r] & H^1(k, \widehat{T})^D \ar[d] & \\
0 & H^0(k,\widehat{T})^D \ar[l] & P^2(k, T) \ar[l] & H^2(k, T) \ar[l] &
}
\end{displaymath}

Et la suite duale correspond exactement \`a la suite de Poitou-Tate pour le module des caract\`eres de $T$ :
\begin{displaymath}
\xymatrix{
0 \ar[r] & \widehat{T}(k)^{\wedge} \ar[r] & P^0(k, \widehat{T})^{\wedge} \ar[r] & H^2(k, T)^D \ar[d] \\
& H^1(k, T)^D \ar[d] & P^1(k, \widehat{T}) \ar[l] & H^1(k, \widehat{T}) \ar[l] \\
& H^2(k, \widehat{T}) \ar[r] & P^2(k, \widehat{T})_{\textup{tors}} \ar[r] & (T(k)^D)_{\textup{tors}} \ar[r] & 0
}
\end{displaymath}
et fournit l'isomorphisme $H^3(k, \widehat{T}) \cong P^3(k, \widehat{T})$.
	\item Si $G$ est semi-simple simplement connexe, dans le cas non compact, cette suite se d\'ecoupe en trois isomorphismes :
$$ \overline{G(k).G_{S_0}} = P^0(k, G)$$
qui n'est autre que le th\'eor\`eme d'approximation forte pour $G$,
$$H^1(k, G) \cong P^1(k, G)$$
qui traduit le principe de Hasse et l'approximation faible pour $G$, et 
$$H^2(k, G) / \sim \cong \bigoplus_v \left( H^2(k_v, G) / \sim \right)$$
qui traduit le principe local-global pour le $H^2$ non ab\'elien.
Quant \`a la suite duale, tous ses termes sont nuls.
	\item Si $G$ est semi-simple (toujours dans le cas non compact), de groupe fondamental $B$, la suite de Poitou-Tate devient
\begin{displaymath}
\xymatrix{
& H^1(k, \widehat{B})^D \ar[d] & P^0(k, G) \ar[l] & \overline{G^{\textup{sc}}_{S_0} G(k)} \ar[l] & 0 \ar[l] \\
& H^1(k, G) \ar[r] & P^1(k, G) \ar[r] & \widehat{B}(k)^D \ar[r] & 0 \\
}
\end{displaymath}
et un isomorphisme $H^2(k, G) / \sim \cong \bigoplus_v \left( H^2(k_v, G) / \sim \right)$.
Quant \`a la suite duale, c'est exactement la suite de Poitou-Tate
pour le module fini $\widehat{B}$ :
\begin{displaymath}
\xymatrix{
& H^2(k, B)^D \ar[d] & P^0(k, \widehat{B}) \ar[l] & \widehat{B}(k) \ar[l]  & 0 \ar[l] \\
& H^1(k, \widehat{B}) \ar[r] & P^1(k, \widehat{B}) \ar[r] & H^1(k, B)^D \ar[d] & \\
0 & B(k)^D \ar[l] & P^2(k, \widehat{B}) \ar[l] & H^2(k, \widehat{B}) \ar[l] & 
}
\end{displaymath}
avec l'isomorphisme $H^3(k, \widehat{B}) \cong P^3(k, \widehat{B})$.
	\item Si $G$ v\'erifie les hypoth\`eses du th\'eor\`eme
          \ref{PTNA}, si $T = G^{\textup{tor}}$ d\'esigne le quotient
          torique de $G$, alors les suites de Poitou-Tate de $G$ et de
          $T$ sont compatibles au morphisme $G \rightarrow T$.
\end{itemize}
\end{exs}

\begin{dem}
On remarque d'abord que l'on peut supposer $G$ r\'eductif pour prouver
le th\'eor\`eme.

L'exactitude des deux premi\`eres lignes r\'esulte de la section
pr\'ec\'edente (voir th\'eor\`eme \ref{theo fin}) et de la suite exacte de
Kottwitz-Borovoi d\'emontr\'ee dans \cite{BorAMS}, th\'eor\`eme 5.16,
et qui se d\'eduit de la suite de Poitou-Tate du
th\'eor\`eme 6.1 de \cite{Dem1} par ab\'elianisation. 

Reste \`a montrer l'exactitude de la suite (d'ensembles point\'es) suivante :
$$0 \rightarrow K \rightarrow H^2(k, G) / \sim \rightarrow \bigoplus_v \left(
H^2(k_v, G) / \sim \right) \rightarrow \H^{-1}(k, \widehat{C})^D \rightarrow 0$$
et \`a identifier le noyau $K$ (qui est seulement un ensemble point\'e
a priori) avec le groupe ab\'elien $ \textup{\cyr{SH}}^2(C) $.

Tout d'abord, montrons que la fl\`eche de localisation $H^2(k, G) / \sim \rightarrow \prod_v H^2(k_v, G) / \sim$ a bien son image contenue dans $\bigoplus_v H^2(k_v, G) / \sim$. Prenons $\eta \in H^2(k, G)$, et notons $\eta' := \textup{ab}^2(\eta) \in H^2_{\textup{ab}}(k,G) = \H^2(k,C)$. Alors $(\eta'_v) \in \P^2(k, C)$, et on sait que $\P^2(k,C) = \bigoplus_v \H^2(k_v,C)$ car $H^i(\mathcal{O}_v, \mathcal{S}) = 0$ pour tout $\mathcal{O}_v$-tore $\mathcal{S}$ et tout $i \geq 2$, car $H^i(\mathcal{O}_v, \mathcal{S}) \cong H^i(\F_v, \mathcal{S} \times_{\mathcal{O}_v} \F_v)$ et $\F_v$ est de dimension cohomologique $1$. Ainsi pour presque toute place $v$, $\textup{ab}^2(\eta_v) = 0 \in H^2(k_v, C)$, donc par le th\'eor\`eme 5.5 de \cite{Bor2}, $\eta_v \in H^2(k_v, G)$ est neutre pour presque toute place $v$, donc $(\eta_v) \in \bigoplus_v H^2(k_v, G) / \sim$. Donc l'application $H^2(k,G) / \sim \rightarrow \bigoplus_v H^2(k_v, G) / \sim$ est bien d\'efinie. 
On remarque ensuite que la fl\`eche d'ab\'elianisation
$\textup{ab}^2$ envoie $K$ dans $\textup{\cyr{SH}}^2(C)$. En outre, par
d\'efinition de la relation $\sim$ sur $H^2(k, G)$, l'application $K
\xrightarrow{\textup{ab}^2} \textup{\cyr{SH}}^2(C)$ est injective (ce qui
est plus fort que de dire qu'elle est de noyau trivial). Montrons sa surjectivit\'e. Soit $\alpha' \in \textup{\cyr{SH}}^2(C)$, et consid\'erons le diagramme suivant :
\begin{eqnarray}
\label{diag centre}
\xymatrix{
H^2(k, Z) \ar[r] \ar[d]^{j_*} & P^2(k, Z) \ar[d] \\
\H^2(k, C) \ar[r] \ar[d] & \P^2(k, C) \ar[d] \\
H^3(k, Z^{\textup{sc}}) \ar[r]^{\simeq} & P^3(k, Z^{\textup{sc}})
}
\end{eqnarray}
Par hypoth\`ese, $\alpha'$ s'envoie sur $0$ dans $\P^2(k, C)$, donc aussi dans $H^3(k, Z^{\textup{sc}})$. Par exactitude de la premi\`ere colonne, $\alpha'$ se rel\`eve donc en $\tilde{\alpha} \in H^2(k, Z)$. Alors $\textup{ab}^2(\tilde{\alpha} + n(G)) = \alpha'$, et si l'on regarde l'image $\beta$ de $\tilde{\alpha} + n(G) \in H^2(k, G)$ dans $P^2(k, G)$, on constate que $\textup{ab}^2(\beta) = 0$ dans $\P^2(k, C)$, et donc $\beta$ est neutre localement partout, ce qui signifie exactement que $\alpha := \tilde{\alpha} + n(G) \in K$. On a donc trouv\'e $\alpha \in K$ relevant $\alpha'$. Donc l'application $\textup{ab}^2 : K \xrightarrow{\simeq} \textup{\cyr{SH}}^2(C)$ est une bijection.
Montrons maintenant l'exactitude en $ \bigoplus_v H^2(k_v, G) /
\sim$. Soit donc $\alpha \in \bigoplus_v H^2(k_v, G) / \sim$ d'image nulle dans $\H^{-1}(k, \widehat{C})^D$. On regarde le diagramme suivant :
\begin{displaymath}
\xymatrix{
H^2(k, G) / \sim \ar[r] \ar[d]  & \bigoplus_v H^2(k_v, G) / \sim \ar[r] \ar[d] &  (k[G]^* / k^*)^D \ar[d]^{=} \\
\H^2(k, C) \ar[r] & \P^2(k, C) \ar[r] & \H^{-1}(k, \widehat{C})^D
}
\end{displaymath}
dont on sait que la seconde ligne est exacte (voir th\'eor\`eme 6.1 de
\cite{Dem1}). Soit alors $(\alpha_v) \in \prod_v H^2(k_v, G)$
relevant $\alpha$, donc d'image nulle dans $ (k[G]^* / k^*)^D$. On
pose $(\alpha'_v) := \textup{ab}^2(\alpha_v) \in \P^2(k, C)$. Par
exactitude de la seconde ligne, $(\alpha'_v)$ provient d'un
\'el\'ement $\alpha'$ de $\H^2(k, C)$. Comme plus haut, on \'ecrit
$(\alpha_v) = (\eta_v) + (n(G)_v)$ avec $(\eta_v) \in P^2(k, Z)$ :
alors avec le diagramme (\ref{diag centre}), $\alpha'$ s'envoie sur
$0$ dans $H^3(k, Z^{\textup{sc}})$, donc $\alpha' = j_*(\eta^0)$ pour un
certain $\eta^0 \in H^2(k, Z)$. Alors l'\'el\'ement $\alpha^0 :=
\eta^0 + n(G) \in H^2(k, G) / \sim$ s'envoie sur $\alpha \in
\bigoplus_v H^2(k_v, G) / \sim $ par injectivit\'e de la fl\`eche
$\bigoplus_v H^2(k_v, G) / \sim \rightarrow \P^2(k, C)$. 

Reste \`a montrer la surjectivit\'e de la fl\`eche $\bigoplus_v
H^2(k_v, G) / \sim \rightarrow (k[G]^* / k^*)^D$. Pour cela, on
utilise la surjectivit\'e de la fl\`eche $\P^2(k, C) \rightarrow
H^{-1}(k, \widehat{C})^D$ (voir th\'eor\`eme 6.1 de \cite{Dem1}). Or pour
toute place $v$, l'image de $\textup{ab}^2 : H^2(k_v, G)
\rightarrow \H^2(k_v, C)$ est exactement $\textup{Im}(H^2(k_v, Z)
\xrightarrow{j_*} \H^2(k_v, C))$; cette image est exactement
$\H^2(k_v, C)$ si $v$ est une place finie, mais elle peut \^etre plus
petite si $v$ est infinie. On conclut de la fa\c con suivante : soit
$\delta \in \H^{-1}(k, \widehat{C})^D$. On sait qu'il existe $c \in
\P^2(k, C)$ d'image $\delta$ dans $\H^{-1}(k, \widehat{C})^D$. Notons
$c_{\infty} \in \prod_{v \in S_{\infty}} \H^2(k_v, C)$ la composante
de $c$ aux places infinies. On sait que le morphisme $\H^2(k, C)
\rightarrow \prod_{v \in S_{\infty}} \H^2(k_v, C)$ est surjective (par
d\'evissage, en se ramenant au cas des tores). Ainsi $c_{\infty}$ se
rel\`eve en $c' \in \H^2(k, C)$. Notons $c'' := c - c' \in \P^2(k, C)$. Alors
$c''_{\infty} = 0$ et $c''$ s'envoie sur $\delta$ dans $\H^{-1}(k,
\widehat{C})$. Alors $c''$ se rel\`eve dans $\bigoplus_v
H^2(k_v, G) / \sim$ en un \'el\'ement $g$ dont l'image dans $(k[G]^* /
k^*)^D$ est exactement $\delta$, ce qui conclut la preuve.
\end{dem}

\section{Points entiers sur les torseurs sous un groupe lin\'eaire}
\label{section entier}
Dans cette section, on applique les r\'esultats de la section
\ref{subsection AF} pour \'etudier l'obstruction de Brauer-Manin
enti\`ere sur un torseur sous un groupe lin\'eaire connexe.
\begin{theo}
\label{theo entier}
Soit $\mathcal{X}$ un $\mathcal{O}$-sch\'ema plat, dont la fibre
g\'en\'erique $X$ est un torseur sous un $k$-groupe lin\'eaire connexe $G$. Soit $S_0$ un ensemble fini de places de $k$ tel que $G^i_{S_0} := \prod_{v \in S_0} G^i(k_v)$  est non compact pour tout $k$-facteur presque $k$-simple $G^i$ de $G^{\textup{sc}}$.

Soit alors un ensemble fini $S$ de places de $k$ contenant $S_0$. On suppose qu'il existe un point $(P_v) \in \mathcal{X}(\A_S)$ qui est orthogonal au groupe $\textup{Br}_1(X)$. Alors il existe un point dans $\mathcal{X}(\mathcal{O}_S)$ qui est arbitrairement proche de $P_v$ pour $v \in S \setminus S_0$ non archim\'edienne, et dans la m\^eme composante connexe de $X(k_v)$ pour $v \in S$ r\'eelle.

En particulier, si $S_0 \subset S_{\infty}$ et si $(P_v) \in \mathcal{X}(\A_{\mathcal{O}})$, alors $\mathcal{X}(\mathcal{O}) \neq \emptyset$.
\end{theo}

\begin{dem}
On sait que $X(\A_k)^{\textup{Br}_a} \neq \emptyset$, donc par le corollaire 2.5 de \cite{Bor}, $X(k) \neq \emptyset$. Donc $X$ est $k$-isomorphe \`a $G$. On conclut alors gr\^ace au corollaire \ref{cor fin}, \`a la mani\`ere du th\'eor\`eme 4 de \cite{HarAF}.
\end{dem}

\addcontentsline{toc}{chapter}{Bibliographie}
\bibliography{AFgroupes}
\bibliographystyle{amsalpha}

\end{document}